\documentclass[12pt]{article}

\usepackage{amssymb,amsmath}

\usepackage{graphicx}

\renewcommand{\a}{\alpha}                 
\renewcommand{\b}{\beta}                  
\newcommand{\g}{\gamma}                   
\newcommand{\de}{\delta}                  
\newcommand{\D}{\Delta}                   
\newcommand{\ep}{\varepsilon}             
\newcommand{\la}{\lambda}                
\newcommand{\si}{\sigma}                  

\newcommand{\bZ}{\mbox{\boldmath$Z$}}     
\newcommand{\bM}{\mbox{\boldmath$M$}}     
\newcommand{\bP}{\mbox{\boldmath$P$}}     
\newcommand{\bT}{\mbox{\boldmath$T$}}     
\newcommand{\bA}{\mbox{\boldmath$A$}}     
\newcommand{\bI}{\mbox{\boldmath$I$}}     
\newcommand{\bD}{\mbox{\boldmath$D$}}     
\newcommand{\bQ}{\mbox{\boldmath$Q$}}     
\newcommand{\bnull}{\mbox{\boldmath$0$}}     
\newcommand{\bun}{\mbox{\boldmath$1$}}     
\newcommand{\bsA}{\mbox{\scriptsize\boldmath$A$}}     
\newcommand{\bu}{\mbox{\boldmath$u$}}     
\newcommand{\br}{\mbox{\boldmath$r$}}     
\newcommand{\bff}{\mbox{\boldmath$f$}}     
\newcommand{\bx}{\mbox{\boldmath$x$}}     
\newcommand{\be}{\mbox{\boldmath$e$}}  
\newcommand{\bv}{\mbox{\boldmath$v$}}     

\newcommand{\exn}{\mathbb{E}}                   

\begin{document}

\title{High host density favors greater virulence: \\
a model of parasite-host dynamics based on multi-type branching processes
}

\author{K.~Borovkov,\footnote{Department of Mathematics and Statistics, The University of Melbourne, Parkville 3010, Australia; e-mail: {borovkov@unimelb.edu.au}}
 \
 R.~Day,\footnote{Department of Zoology, The University of Melbourne, Parkville 3010, Australia; e-mail:
 {r.day@unimelb.edu.au}}
 \
 T.~Rice\footnote{Department of Mathematics and Statistics, The University of Melbourne, Parkville 3010, Australia; e-mail:
 {t.rice@ms.unimelb.edu.au}}}

\date{}

\maketitle

\begin{abstract}
We use a multitype continuous time Markov branching process model to describe the
dynamics of the spread of parasites of two types  that can mutate into each other in a
common host population. Instead of using a single virulence characteristic which is
typical of most mathematical models for infectious diseases, our model uses a
combination of two characteristics: lethality and transmissibility. This makes the
model capable of reproducing the empirically observed fact that the increase in the
host density can lead to the prevalence of the more virulent pathogen type. We provide
some numerical illustrations and discuss the effects of the size of the enclosure
containing the host population on the encounter rate in our model that plays the key
role in determining what pathogen type will eventually prevail. We also present a
multistage extension of the model to situations where there are several populations and
parasites can be transmitted from one of them to another.
\end{abstract}

\section{Introduction}
\label{Sec_1}

There is an extensive literature on mathematical modeling of infectious diseases. Both
the rate of spread and the virulence of pathogens are important, and both the dynamics
of disease spread and the evolution of virulence in parasite-host systems have
attracted much attention. The first mathematical model of disease dynamics apparently
goes back to the 1760 D.~Bernoulli's paper~\cite{Be} (see also pp.\,2--6
in~\cite{DaGa}). A major step in further development of deterministic modelling was the
``threshold theorem"~\cite{KeMc} showing that the initial frequency of susceptibles
must exceed a certain threshold value to give rise to an epidemic. The first stochastic
epidemiological models appeared in the late 1920s~\cite{Mc}, applying the ``law of mass
action"  suggested in~\cite{Be} to probabilistic description of the epidemics.
Stochastic models became much more popular after Bartlett~\cite{Bar0} formulated a
model for the general stochastic epidemic by analogy with the deterministic model
from~\cite{KeMc}. May and Anderson~\cite{MaAn} 
introduced the first model of evolutionary change in pathogen virulence in 1983, which
explained the decline in virulence of myxomatosis, introduced to control Australian
rabbits, in terms of optimal transmission of the parasite between hosts. Subsequent
work, notably~\cite{MaAn, Ew1, Ma1, Kn, GaEtal}, 
has expanded
on this classic paper, to cover the effects of a vector of the pathogen, vaccination of
hosts, and various other issues. For more extensive surveys of the field and literature
reviews, we refer the reader to the monographic literature~\cite{AnMa1, AnBr, Bai, Bec,
DaGa, MoSl} and also to papers~\cite{Al, DiSc, He}.

In natural systems, the density of host organisms declines when a virulent
pathogen produces an epidemic that kills host rapidly. The models of virulence
have made clear that this leads to selection for less virulent strains, such
that infected hosts live longer, and thus tend to pass the pathogen to more
susceptible new hosts, maintaining the epidemic.  While hosts will have longer
generation times than their pathogens, they too, will evolve to be more
resistant, further reducing the virulence of the pathogen for these hosts. Thus
one can expect pathogens to be highly virulent only when they have recently
switched hosts, such as the SARS coronavirus that was probably transferred to
humans from animals such as Himalayan civet cats~\cite{GuEtAl} or when the
plague bacterium, {\em Yersina pestis}, has been transferred to humans from
rodents~\cite{Ye}.

The main objective of the present paper is to provide quantitative and also qualitative
(cf.\ the discussion of the above-mentioned threshold theorem and other findings
from~\cite{KeMc} on pp.\,11 and~29 in~\cite{DaGa}) insights into the dynamics of
pathogen populations where the virulence of the pathogens can change due to mutation.
More specifically, we are interested in analyzing the effects that {\em changes in host
density\/} may have on the virulence of the pathogen. In agricultural systems, humans
have been maintaining animals and plants at very high densities for about 10,000 years,
and humans themselves began living at high local densities in villages, towns and
cities from about the same time~\cite{Cr}. Densities of both humans and domesticated
organisms have probably increased continually since then, but the densities at which
some domestic animals are kept appear to have increased very rapidly in the last few
decades (see e.g.~\cite{Fra}, and the recent rapid development of large scale
commercial aquaculture involves the maintenance of very high densities of a whole new
set of marine species that will have their own pathogens
(see~\cite{HaEtal1,HaEtal,Beh}).  Thus the question arises as to whether high host
densities are likely to lead to selection for more virulent strains of pathogens.


The mathematical models we discuss in the present paper incorporate a continuous time
Markovian multitype branching process. Branching processes first appeared in
epidemiological context in the mid-1950s: one can mention here the paper~\cite{Bh} and
also refer to ~\cite{Bar, Ke, Wh}. For further literature review, we again refer the
reader to the above-listed monographic and survey literature, while some  examples of
applications of multitype branching processes in epidemiology can be found
in~\cite{BecMa, Mee, Hei}.

Branching processes are relatively simple and well-understood mathematical models; for
treatments of the theory of branching processes, see e.g.~\cite{Ja, AtNe}, a recent
addition to the monographic literature on branching processes prepared specifically for
biological audience being~\cite{HaJaVa}. At the same time, the processes were shown to
be good approximations to the general stochastic epidemic models at the initial stage
(when the total population size is large and the initial number of infectives is small)
and also at the final stages of the epidemics (see e.g.~\cite{Lu, BaDo, He}). As our
main objective is not to model the detailed dynamics and describe the whole history of
the epidemic itself, but rather to suggest a model for changes in the pathogen's
population composition depending on the   density of the host population, we will
restrict ourselves to working with the simpler branching process models.

Moreover, the most critical stage of an epidemic is the initial one, when it is
basically determined if there will be a large-scale event or the epidemic will die out.
And as branching processes are good approximations to the general stochastic epidemic
models at the initial stage, the threshold analysis aimed to determine if the  ``basic
reproductive number" (defined roughly as the expected number of secondary cases
produced by one infected individual) is greater than one (which implies the danger of
an outburst or persistence of endemic levels) can be carried out using those models
(see e.g.\ Chapters~6 and~8 in~\cite{MoSl}).

An additional argument for this approach is that in agriculture or aquaculture
situations, when a threat of an epidemic outburst arises, drastic measures are
usually taken: the affected part of the population may be isolated or even
destroyed. These measures will substantially change the way pathogens can pass
between hosts, and hence can make the standard stochastic epidemic models
inapplicable beyond the initial stage of the epidemic anyway.

One can also envisage an extension of our simple analysis to the general stochastic
models as the same mechanism will certainly work for the latter as well. However, such
extensions will be much less tractable analytically and may lead to no closed-form
answers.

In our analysis, we will look at   supercritical two-type branching processes (so that
the basic reproductive number will be greater than one: we are interested in what
happens when {\em there is\/} a threat of epidemics) and then look at the behaviour of
the ratio of the sizes of the subpopulations in the process (representing two versions
of the pathogen, that can mutate into each other). This quantity can be used to
determine which of the two types will become dominant in the population over time. One
of the advantages of the branching process model is that one can easily incorporate in
it two different characteristics of pathogenes' virulence: their {\em lethality\/}
(which can be described by the mean survival time of an infected host) and {\em
transmissibility\/} (specifying the probability of an infected host infecting a
susceptible one on their contact), instead of representing them together by a single
quantity termed simply ``virulence" (cf.\ e.g.~\cite{AnMa}). For recent discussions of
the interplay between these two characteristics, see~\cite{GiEtal,LiMo}.

In the paper, we use our approach to model the dynamics of a host/parasite
population where parasites can be of one of two types that can differ in their
lethality and transmissibility. The underlying simple continuous time Markov
model of a two-type branching process is presented in Section~\ref{Sec_2}. The
analysis of the model and derivation of the dynamics for the mean functions are
given in Section~\ref{Sec_3} and used in Section~\ref{Sec_4} to establish the
eventual composition of the pathogen population. It turns out that our model is
capable of reproducing the above-mentioned phenomenon of shifting the overall
parasite's lethality in response to increased density of the hosts.  This is
actually  possible because of using two different parameters to characterise
the virulence of the pathogen.

In Section~\ref{Sec_5} we present a few remarks on how the change in the size
of the enclosure a given population of hosts inhabits can affect the
``encounter rate" for the hosts\,---\,the key parameter of the model describing
the ``effective density" in the host population. Finally, in
Section~\ref{Sec_6} we consider a multi-stage modification of our model that
can be used to analyse farm or aquaculture situations in which there are many
enclosures or tanks of animals, and an outbreak of an infectious disease occurs
in one of them. We assume that the pathogen might at some stage be transmitted
to the next enclosure, where the epidemic process starts anew etc., and discuss
possible scenarios of the development of the epidemic in the farm.
Section~\ref{Sec_7} contains a few final remarks concerning biological
interpretation of our results.

\section{Description of the branching process model}
\label{Sec_2}

Assume that we have a large population of hosts that can be infected by
parasites of one of two types that will be denoted by $T_1$ and $T_2$. The
pathogen types can differ in both their lethality   and transmission rate. The
numbers of infected hosts at time $t$ are represented by the vector
\[
\bZ (t) = (Z_1 (t), Z_2 (t)), \quad t\ge 0,
\]
$Z_i(t)$ being the time $t$ number of hosts infected with the type  $T_i$ pathogen (for
simplicity,  at any given time, any given host can be infected with one type of the
pathogen only). We do not keep track of the number of hosts that remain uninfected
(susceptible), assuming instead that this number will remain large enough during the
time period for which our mathematical model is intended, so that the dynamics of the
process $\{\bZ (t),\, t\ge 0\}$ to be described do not change over time.

A general description of the model is as follows. We assume that each infected
individual lives a random time (which will tend to be shorter when one is
infected with the ``more lethal" of the two pathogen types). During its
lifetime, an infected host can encounter susceptible hosts and, with a
probability depending on the type of the pathogen it carries, transmit the
parasite to them. The rate of such (random) encounters will be specified by a
special parameter that we can vary in order to model changes in the density of
the host population.

Finally, we also allow the pathogens to mutate, so that when a host originally infected
with $T_1$ encounters a susceptible host, the latter can become infected with
$T_2$-type parasites (and the other way around).

Now we will present a formal mathematical model. First recall that the exponential
distribution with rate (or intensity) $\alpha >0$ has density of the form
\begin{equation}
p(t)= \alpha e^{-\alpha t},\qquad t>0,
 \label{Exp_den}
\end{equation}
with mean  $1/\alpha$, and plays a special role in probability theory due to
its unique memoryless property that makes the distribution ubiquitous in the
theory of Markov processes. Namely, a random variable $\tau>0$ modelling, say,
the time of the first encounter of a given infected host with a healthy one,
has this property if, for any $s,t>0$, the conditional probability of the event
$\{\tau>s+t\}$ given that $\tau>s$ coincides with the probability of
$\{\tau>t\}$:
\begin{equation}
\Pr ( \tau>s+t \,|\,  \tau>s )
 \equiv \frac{\Pr ( \tau>s+t ,\,  \tau>s )}{\Pr (  \tau>s )}
 \equiv \frac{\Pr ( \tau>s+t )}{\Pr (  \tau>s )}
 = \Pr (\tau
>t).
 \label{Mem_less}
\end{equation}
In words, if, at time $s$ we know that there has been no such encounter, then the
conditional distribution (given that information) of the residual random time $\tau-s$
till the encounter will be the same as the original distribution of $\tau$. It is
obvious that if $\tau$ has density~\eqref{Exp_den} then
\begin{equation}
\Pr (\tau >t ) = \int_t^\infty p(s) ds = e^{-\alpha t},\qquad t>0,
 \label{Exp_tail}
\end{equation}
and so the property \eqref{Mem_less} is clearly satisfied.

An equivalent formulation of the property can be given in terms of the
distribution's hazard rate $r_\tau(s)$ that quantifies the probability that,
given there has been no encounter up to time $s$, there will  be one
``immediately afterwards": in case a random variable $\tau$ has a continuous
density $p_\tau,$ the hazard rate is defined by
\begin{equation}
r_\tau (s) =\frac{p_\tau (s)}{\Pr (  \tau>s )}
  \equiv - \frac{d}{ds} \ln \Pr (  \tau>s ) , \qquad s>0,
 \label{Haz_rate}
\end{equation}
and, as $h\downarrow 0,$
\begin{equation}
\Pr (\tau \le s + h |\, \tau >s) = r_\tau (s) h + o(h), \qquad s>0,
 \label{Haz_prob}
\end{equation}
where, as usual, $o(h)$ denotes a quantity that vanishes faster than $h$: $o(h)/h\to
0.$

It is obvious from~\eqref{Exp_tail} and \eqref{Haz_rate} that the hazard rate of a
distribution is constant if and only if it is exponential (in that case, the hazard
rate simply equals the distribution's rate). In applications, one often uses
exponentially distributed random variables to model times between successive events of
a particular kind and also lifetimes. This is because, on the one hand, such
assumptions make sense from the modelling view point (in a large population, meeting a
new individual during a time interval $[t,t+h]$, $h>0$, can scarcely depend on one's
``life history" prior to time $t$) and, on the other hand, as the resulting models are
usually Markovian, so that the powerful machinery of the theory of Markov processes is
applicable.

Our basic model assumptions are as follows:
\begin{itemize}
\item[(a)]
The initial population contains $z_i$ hosts infected with parasites of type $T_i$,
$i=1,2.$

\item[(b)] A host infected with type $T_i$ parasites lives a random time
    which is exponentially distributed with parameter $\alpha_i >0,$
    $i=1,2.$ Clearly, the pathogen with a higher rate $\alpha_i$ will be
    the more lethal one, as the mean lifetime in that case will be lower.

\item[(c)] Any infected host can encounter susceptible ones. The time till
    the first encounter of a given infected host (of any type) with a
    susceptible host is a random variable exponentially distributed with
    rate $\la$. It is clear from the above discussion of the properties of
    exponential distributions that, at any time $t$, the residual times
    till encounters of the current infected hosts with susceptible hosts
    are all exponential with the same rate~$\la$. A similar observation
    applies to the residual lifetimes; this ensures that the process $\{\bZ
    (t)\}$ will be {\em Markovian\/}: given the current state of the
    process, its future evolution does not depend on the past one.

In a modification of the model, one can assume that, at time $t$, the time
till the first encounter of a given infected host (of any type) with a
susceptible host has a hazard rate $\la (t)$ which {\em can depend
on\/}~$t$. This will enable one to model changes in the density of the host
population that occur over time (the higher $\lambda$, the more often are
encounters, which corresponds to higher host density situations). The
process will remain Markovian, but will become time-inhomogeneous.

\item[(d)] At any encounter with susceptible hosts, a  $T_i$-infected host
    meets only one susceptible host (there can be several such encounters
    during the host's lifetime). At each such instance, the $T_i$-infected host
    transmits the parasites to the susceptible one with probability
    $\beta_i$ (so that, with the complement probability
    $1-\beta_i$, the encounter will have no consequences for the
    susceptible host).

\item[(e)]
Mutations $T_1 \leftrightarrow T_2$ are possible. A host infected with $T_1$-type
pathogens will remain such till its end, but, when transmitting pathogens to a
susceptible host during an encounter, the newly infected host will carry $T_2$-type
pathogens with a probability $\mu_1$. Likewise, $\mu_2$ denotes the probability that a
successful transmission of parasites from a $T_2$-infected host to a susceptible one
resulted in making the latter $T_1$-infected.

\item[(f)]
All the above-mentioned random times (lifetimes, times till the first encounter) are
independent of each other.
\end{itemize}

Of course, the above assumptions oversimplify the real biological processes. There are
likely to be several or many strains of pathogen, and the probability of mutation from
one to another will vary. We suggest that simplifying the system to two strains is
likely to retain the same key dynamic features. Further, the distributions of times
until events occur are likely to be approximately exponential as argued earlier, and we
do not see any reasons why host survival times and encounter rates should depend in any
way on each other. Note that ``encounters" between hosts are simply occasions where a
pathogen can be transferred between hosts. They do not have to involve physical
contact. Thus a pathogen transferred by aerosols might be transferred between pigs that
are isolated in neighbouring stalls.

\begin{figure}[ht]
 \label{Fig_1}
\setlength{\unitlength}{1 mm}
\begin{center}
\begin{picture}(75, 28)(0,2)
 \put(5,5){\vector(1,0){65}}
 \put(5,5){\line(0,1){25}}
%
%

 \put(25,4){\line(0,1){2}}
 \put(53,4){\line(0,1){2}}
{\small 
 \put(4.5,1){$0$}
 \put(23,0.5){$1/\la'$}
 \put(51,0.5){$1/\la''$}
 \put(69,0.5){time}
%
%
 {\thicklines \put(5,15){\line(1,0){25}}
  \put(5,25){\line(1,0){45}}}
   \put(30,14){\line(0,1){2}}
   \put(50,24){\line(0,1){2}}
 \put(28,10){$1/\a_2$}
 \put(48,20){$1/\a_1$}
%
%
}
\end{picture}
\end{center}
\caption{Mean lifetimes of infected hosts and mean times to encounter.}
\end{figure}
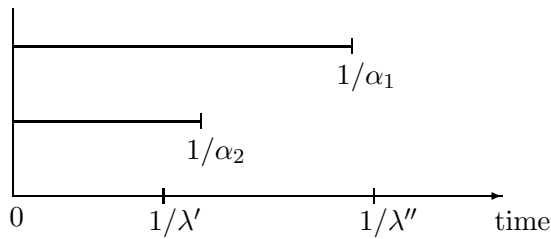

The diagram in Fig.~1 illustrates the ``physical" meaning of our assumptions. The two
horizontal segments represent the mean lifetimes of hosts infected by the pathogens of
our two types: the longer segment (of length $1/\a_1$) corresponds to $T_1$ which we
assume less lethal by stipulating that $\a_2 >\a_1$, whereas the shorter one (of length
$1/\a_2$) corresponds to $T_2$.

When the host population density is relatively low (say, represented by the
value $\la =\la''$, as depicted), the lifetime of $T_2$-infected hosts will be
too short compared to the mean time $1/\la''$ between encounters to give them a
good opportunity to encounter susceptibles and hence further propagate in the
host population. One may expect that this will result in the eventual
prevalence of $T_1$ pathogens who have better chance of being transmitted as
they live longer. However, if the host population density increases (say, to
$\la =\la'$, as depicted), then the more lethal type $T_2$ may have frequent
enough encounters which, combined with its higher transmissibility, can lead to
its eventual prevalence. As we will see in Section~\ref{Sec_3}, the above
argument is confirmed by mathematical analysis.


Assumptions (a)--(f) imply that our process $\{\bZ (t)\}$ is actually a {\em
two-type time homogeneous Markov branching process in continuous time\/}, see
e.g.\ Section~V.7 in~\cite{AtNe}. That is, $\{\bZ (t)\}$ describes the dynamics
of a population consisting of individuals (or ``particles") of two types, 1
and~2. The transitions of different particles in the process are assumed to be
independent. A particle of type $i$ lives for an exponentially random time with
rate $a_i$. At the end of its life it disappears. It can either
\begin{itemize}
\item[(i)] simply disappear (in terms of our modelling assumptions above, this means that
a given $T_i$ infected host died having never encountered a susceptible host), or

\item[(ii)] produce one particle of the same type (meaning: there was an encounter, but no
transmission occurred; we think of the ``newly produced" particle of type $i$ as just
the ``old" infected host of type~$T_i$ who keeps living --- note that, due to the
memoryless property of the exponential distribution, such an identification of a ``new"
particle with the ``old" one is in agreement with our assumptions (a)--(f) (so that the
life of one $T_i$-infected host can actually be represented by a succession of several
type $i$ particles, for which reason we do not use $T_i$ to denote the type of
particles in the branching process), or

\item[(iii)] produce two particles of the same type $i$ (meaning: there was an encounter
and successful transmission, but no mutation; one of the ``newly produced" particles is
actually the original host, the other represents a newly infected --- with the same
$T_i$-type pathogen --- host), or

\item[(iv)] produce two particles of different types (meaning: there was an
    encounter and successful transmission {\em and\/} mutation; one of the
    ``newly produced" particles is actually the original host, the other
    represents a newly infected --- with the pathogen of the other type ---
    host).
\end{itemize}

To see that both sets of assumptions (a)--(f) and (i)--(iv) describe the same
dynamics of $\{\bZ(t)\}$, it suffices to note that in both cases we deal with
time-homogeneous continuous time jump Markov processes (which follows from the
exponentiality assumptions) and then verify that, choosing suitable parameter
values for the second model, one can obtain the same transition rates as for
the first one.

To do that, we first observe that assumptions (i)--(iv) imply the {\em branching
property\/} which means that, for any $s\ge 0$, given $\bZ(s) = (z_1, z_2),$ the {\em
future\/} $\{\bZ (s+t),\, t>0\}$ of the process will follow the same probability laws as
that of the sum of $z_1$ independent copies of $\bZ (t)$ starting at time~0 with a
single particle of type~1 and $z_2$ independent copies of $\bZ (t)$ starting at time~0
with a single particle of type~2. It is clear that the first model (specified by
(a)--(f)) has the same property.

Moreover, the branching property implies that to completely describe the evolution of
the process, it suffices to specify transition probabilities
\[
p_i^{(h)} (j_1, j_2) = \Pr \bigl(\bZ(h) = (j_1, j_2) |\, \bZ(0)=\be_i\bigr)
\]
from the basic initial states
\[
\be_1 = (1,0)\qquad \mbox{and} \qquad  \be_2 = (0,1)
\]
for arbitrary small time increments~$h$. It is obvious that the probability of
having more than one transition during a small time interval $(0,h)$  will
be~$o(h)$, so we just need to consider {\em  where\/} a single transition can
take the process from a basic state $\be_i$ according to assumptions~(a)--(f)
and show that the transitions will have the same rates as for a process
specified by~(i)--(iv) (for a suitable choice of parameter values).

Suppose that, in our branching process, particles of type $i$ have
exponentially distributed random lifetimes with rates $a_i=\a_i+\la$, $i=1,2$.
Moreover, at the end of a  type $i$ particle's life, it produces a random
number of children (possibly of both types) according to the offspring
distributions $q_i (j_1, j_2)=\Pr \bigl($a particle of type $i$ gives birth to
$j_1$ particles of type $1$ and $j_2$ ones of type~2$\bigr)$ given by the
following table:
\[
\begin{array}{c}
\begin{array}{c|cc}
(j_1, j_2)  & q_1 (j_1, j_2) & q_2 (j_1, j_2)\vphantom{\displaystyle\sum} \\
 \hline
  (0,0) &  \displaystyle \frac{\a_1}{\la + \a_1}
 & \displaystyle \frac{\a_2}{\la + \a_2} \vphantom{\int_1}\\
 \hline
  (1,0) & \displaystyle \frac{(1-\b_1)\la}{\la + \a_1}
 &  \displaystyle 0 \vphantom{\int^{A^A}_A}\\
  \hline
  (2,0) & \displaystyle \frac{(1-\mu_1)\b_1 \la}{\la + \a_1}
 &  \displaystyle 0 \vphantom{\int^{A^A}_A}\\
  \hline
  (1,1) & \displaystyle \frac{ \mu_1 \b_1 \la}{\la + \a_1}
 &  \displaystyle\frac{ \mu_2 \b_2 \la}{\la + \a_2} \vphantom{\int^{A^A}_A}\\
  \hline
  (0,2) &  \displaystyle 0
 & \displaystyle \frac{(1-\mu_2)\b_2 \la}{\la + \a_2}\vphantom{\int^{A^A}_A}\\
 \hline
  (0,1) & \displaystyle 0
 &  \displaystyle\frac{(1-\b_2)\la}{\la + \a_2} \vphantom{\int^{A^A}}\\
 \end{array}
  \\
  \\
  \mbox{Table 1. Offspring distributions in the branching process model (i)--(iv).}
   \end{array}
 \]

Then, due to independence and \eqref{Haz_prob},
\begin{align*}
p_1^{(h)} (0,0) & = \Pr \bigl(\mbox{initial type~$1$ particle dies in $(0,h)$, no
children produced}\bigr)
 \\
 & = \Pr \bigl(\mbox{initial type~$1$ particle dies in $(0,h)$}\bigr) \times q_1 (0,0)\\
 & = ((\a_1 + \la)h + o(h)) \times \frac{\a_1}{\la + \a_1} =  \a_1 h + o(h),
 \end{align*}
which is clearly the same as the probability for a single $T_1$-infected host (from the
first model) to die during $(0,h)$.

Likewise,
\begin{align*}
p_1^{(h)} (2,0) & = \Pr \bigl(\mbox{initial type~1 particle dies in $(0,h)$, producing 2
children of type~1}\bigr)
 \\
 & = \Pr \bigl(\mbox{initial type~1 particle dies in $(0,h)$}\bigr) \times q_1 (2,0)\\
 & = ((\a_1 + \la)h + o(h)) \times \frac{(1-\mu_1)\b_1 \la}{\la + \a_1}
  =  (1-\mu_1)\b_1 \la h + o(h).
 \end{align*}
The corresponding transition in the first model is as follows: a $T_1$-infected
host did not die during $(0,h)$, but met a susceptible; the pathogen was
transmitted, no mutation occurred. Due to independence, the probability of this
will be
\[
(1+O(h)) \times (\la h + o(h)) \times \b_1 \times (1-\mu_1)
 = (1-\mu_1)\b_1 \la h + o(h).
\]
Virtually the same argument shows that
\begin{align*}
p_1^{(h)} (1,1) & = \Pr \bigl(\mbox{initial type~1 particle dies in $(0,h)$, producing
one child of each type}\bigr)
 \\
 & = \Pr \bigl(\mbox{initial type~1 particle dies in $(0,h)$}\bigr) \times q_1 (1,1)\\
 & = ((\a_1 + \la)h + o(h)) \times \frac{ \mu_1 \b_1 \la}{\la + \a_1}
  =   \mu_1 \b_1 \la h + o(h),
 \end{align*}
coincides with the probability for a $T_1$-infected host to meet a susceptible
and transmit the mutated form of the pathogen.

Finally, given $\bZ(0)=\be_1$, the most likely state for $\bZ (h)$ is $\be_1$. In our
branching process this occurs when either the original particle lives through the time
interval (with probability $1- (\a_1 + \la)h + o(h)$) or it dies prior to~$h$ producing
a single type~1 particle (of which the probability is $((\a_1 + \la)h + o(h))\times
(1-\b_1)\la/(\la + \a_1) = (1-\b_1)\la h + o(h)$), so that
\[
p_1^{(h)} (1,0) =  1- (\a_1 + \la)h + (1-\b_1)\la h + o(h)
  =  1- (\a_1 + \b_1\la )h + o(h).
\]
In the first model, to get to this state, the initial $T_1$-infected host
survives for $h$ time units and either has no encounters with susceptibles or
has one, but without transmission of a pathogen. The probability of this is
\begin{align*}
(1-\a_1 h &+  o(h))\times [(1-\la h +  o(h)) +  ((1-\b_1)\la h +  o(h))]\\
 &= (1-\a_1 h +  o(h))\times  (1- \b_1 \la h +  o(h)) = 1- (\a_1 + \b_1\la )h + o(h),
 \end{align*}
the same value as for the branching process model.

Of course, we could also work out the probability for $(1,0)$ just as the complement
probability
\[
p_1^{(h)} (1,0) = 1- \bigl(p_1^{(h)} (0,0) + p_1^{(h)} (2,0) + p_1^{(h)} (1,1)\bigr) +
o(h),
 \]
but the presented argument demonstrates the difference between the interpretation of
the elements of the two models ($T_i$-infected hosts in the first model and type~$i$
particles in the second one are not the same, as we noted earlier).

The same calculations are applicable in the case of the initial state $\be_2$, which
leads us to the transition probabilities presented in Table~2 (for $h\to 0$, the
additive terms $o(h)$ being omitted for brevity).
\[
\begin{array}{c}
\begin{array}{c|cc}
(j_1, j_2)  & p^{(h)}_1 (j_1, j_2) & p^{(h)}_2 (j_1, j_2)\vphantom{\displaystyle\sum} \\
 \hline
  (0,0) &    \a_1 h
 &   \a_2 h \vphantom{\int_1^A}\\
 \hline
  (1,0) & 1- (\a_1 + \b_1\la )h
 &   0 \vphantom{\int^{A}_A}\\
  \hline
  (2,0) & (1-\mu_1)\b_1 \la h
 &    0 \vphantom{\int^{A}_A}\\
  \hline
  (1,1) &  \mu_1 \b_1 \la h
 &   \mu_2 \b_2 \la h  \vphantom{\int^{A}_A}\\
  \hline
  (0,2) &   0
 & (1-\mu_2)\b_2 \la h \vphantom{\int^{A}_A}\\
 \hline
  (0,1) &   0
 &  1- (\a_2 + \b_2\la )h \vphantom{\int^{A}}
 \end{array}\\
  \\
  \mbox{Table~2. Transition probabilities for small $h$ values ($o(h)$ terms omitted).}
   \end{array}
  \]

\section{The dynamics of the means}
\label{Sec_3}

Consider the matrix of mean values $\bM (t) = (M_{ij}(t))$, where
\[
M_{ij}(t) = \exn (Z_j(t) | \, \bZ (0) = \be_i), \qquad i,j=1,2,
\]
is the expected number of type~$j$ particles present in the process at time~$t$ given
that the process started at time~$0$ with a single particle of type~$i$. Clearly,
\[
\bM (0) = \bI \equiv \left[
 \begin{array}{cc}
 1 & 0\\
 0 & 1
 \end{array}
 \right],
\]
the identity matrix, and, using the branching and Markov properties of $\{\bZ
(t)\},$ it is easy to demonstrate that $\{\bM (t), \, t\ge 0\}$ possesses the
operator semigroup property:
\begin{equation}
\bM (s+t ) = \bM ( s  )\bM ( t) \qquad\mbox{for any \quad $s,t\ge 0.$}
 \label{Semi}
\end{equation}
Indeed, given that $\bZ(0) = (z_1, z_2),$ the value of $\bZ (s)$ is just the sum of
$z_1$ independent copies of $\bZ (s)$ starting at time~0 with a single particle of
type~1 and $z_2$ independent copies of $\bZ (s)$ starting at time~0 with a single
particle of type~2. As the process is time-homogeneous, we infer that
\begin{equation}
\exn \bigl(\bZ (s+t) | \, \bZ (s) = (z_1, z_2)\bigr)
 = \sum_{i=1}^2 z_i \exn (\bZ (t) | \, \bZ(0) = \be_i)
 = (z_1,z_2) \bM(t).
 \label{Means}
\end{equation}
From here, using the Markov property and the double expectation law for conditional
expectations, we have
\begin{align*}
 \exn \bigl(\bZ (s+t) | \, \bZ(0) = \be_i\bigr)
 &= \exn \bigl[ \exn (\bZ (s+t) | \, \bZ(s))\, \big|\, \bZ(0) = \be_i\bigr]\\
 &= \exn \bigl[  \bZ (s) \bM(t)  \big|\, \bZ(0) = \be_i\bigr] =(M_{i1}(s), M_{i2} (s))\bM(t),
\end{align*}
which is equivalent to~\eqref{Semi}.

Relation~\eqref{Semi} implies (cf.\ p.\,202 in~\cite{AtNe}) that one has the matrix
exponential representation
\begin{equation}
\bM (t) = e^{t\bsA} \equiv \sum_{k=0}^\infty \frac{t^k \bA^k}{k!},\qquad t\ge 0,
 \label{Expl}
\end{equation}
where $\bA^0=\bI$ and
\begin{equation}
\bA = \lim_{h\downarrow 0} \frac{1}{h} (\bM(t) - \bI)
 \label{Gene}
\end{equation}
is the so-called infinitesimal generator of $\{\bM (t)\}$. Indeed,
from~\eqref{Semi} one has $\bM (t+h) - \bM (t) =  (\bM(h)-\bI)\bM ( t)$, so
that~\eqref{Gene} implies that $\dfrac{d}{dt}\bM (t)=\bA \bM (t)$, $\bM
(0)=\bI,$ for which~\eqref{Expl} is clearly a solution, as seen from its
term-wise differentiation.

Evaluating matrix exponentials is rather straightforward: it basically reduces to
calculating the values of the function on the matrix' spectrum (for more detail on
functions of matrices, see e.g.\ Chaper~V in~\cite{Ga}). If, say, $\bA$ is
diagonalisable, so that there exists an invertible matrix $\bQ$ (with inverse
$\bQ^{-1}$: $\bQ^{-1}\bQ = \bQ \bQ^{-1} = \bI$) such that
\begin{equation}
\bQ^{-1} \bA \bQ = \bD \equiv \mbox{diag}\, \{\si_+, \si_-\}
 \equiv \left[
\begin{array}{cc}
 \si_+ & 0\\
 0 & \si_-
 \end{array}
 \right]
 \label{Tran}
\end{equation}
for some $\si_\pm$ (which is the case in our situation, as we will see below), then
clearly $\bA = \bQ \bD \bQ^{-1}$,
$$
\bA^2 = (\bQ \bD \bQ^{-1})^2 = \bQ \bD \bQ^{-1} \bQ \bD \bQ^{-1} = \bQ \bD^2 \bQ^{-1}
 = \bQ \, \mbox{diag}\, \{\si_+^2, \si_-^2\} \, \bQ^{-1}
$$
and so on: $\bA^k = \bQ \bD^k \bQ^{-1}= \bQ \, \mbox{diag}\, \{\si_+^k, \si_-^k\} \,
\bQ^{-1}$. We obtain
\begin{equation}
e^{t\bsA} = \sum_{k=0}^\infty \frac{t^k  (\bQ \bD \bQ^{-1})^k}{k!}
 = \bQ \sum_{k=0}^\infty \frac{t^k \, \mbox{diag}\, \{\si_+^k, \si_-^k\} }{k!} \, \bQ^{-1}
 = \bQ \,  \mbox{diag}\, \{e^{t\si_+},e^{t \si_-}\}  \, \bQ^{-1}.
 \label{Expo}
\end{equation}

Thus, to derive the dynamics of the means, we just have to compute the generator~$\bA$,
which can easily be done using Table~2. Indeed, we infer from the table that, as $h\to
0$,
\begin{align*}
 M_{11} (h) & = 0\times \a_1 h
 + 1 \times [(1-(\a_1+ \b_1 \la) h) + \mu_1 \b_1 \la h]
 + 2 \times (1-\mu_1 )\b_1 \la h
  + o(h)\\
  & = 1 + ((1-\mu_1 )\b_1 \la - \a_1)h + o(h);\\
 M_{12} (h) & = 0\times \a_1 h
 + 1 \times   \mu_1 \b_1 \la h + o(h) = \mu_1 \b_1 \la h + o(h)
\end{align*}
and similarly for $M_{2i} (h)$. From here and~\eqref{Gene} we immediately obtain
\[
\bA = \left[
\begin{array}{cc}
 \g_1 & \de_1\\
 \de_2   &  \g_2
 \end{array}
 \right], \qquad \g_k = (1-\mu_k)\b_k \la - \a_k,
  \quad \de_k =\mu_k \b_k \la,
 \quad k=1,2.
\]
Now solving the characteristic equation $\det(\bA -\si \bI ) = 0$ for $\si$ we find the
eigenvalues $\si_\pm$ of $\bA$ given by
\[
\si_\pm = \frac12 \bigl( \g_1 + \g_2  \pm \D \bigr), \qquad
  \D =\sqrt{ (\g_1 - \g_2)^2
  + 4   \de_1 \de_2   }
\]
(cf.\ similar calculations of the threshold parameter for a somewhat different two-type
model in Section~8.4 of~\cite{MoSl}), with the respective (right) eigenvectors
\[
\bu_\pm = \left[
\begin{array}{c}
 u_\pm\\
 1
 \end{array}
 \right],
 \qquad u_\pm = \frac{\g_1 - \g_2 \pm \D}{2 \de_2} \gtrless 0.
\]
The eigenvalues $\si_\pm$ are clearly different from each other (in any case, this is
guaranteed by the fact that $\si_+$ is the Perron-Frobenius root for the quasi-positive
matrix $\bA$, cf.\ Section~A.8 in~\cite{Th}), which ensures that $\bA$ is diagonalizable
and we can take the transformation matrix $\bQ$ from~\eqref{Tran} to be given by
\[
\bQ  = (\bu_+, \bu_-)
  \equiv
  \left[
\begin{array}{cc}
 u_+ &  u_-\\
 1   &  1
 \end{array}
 \right],
  \qquad
\bQ^{-1} =   \frac{1}{u_+ - u_-}
 \left[
\begin{array}{cc}
 1 &  -u_-\\
 -1  & u_+
 \end{array}
 \right].
\]
Hence \eqref{Expl} and \eqref{Expo} imply that
\begin{align}
 \bM (t) & =  \frac{1}{u_+ - u_-}
  \left[
\begin{array}{cc}
 u_+ &  u_-\\
 1   &  1
 \end{array}
 \right]
   \left[
\begin{array}{cc}
 e^{\si_+ t} &  0\\
 0  &  e^{\si_- t}
 \end{array}
 \right]
 \left[
\begin{array}{cc}
 1 &  -u_-\\
 -1  & u_+
 \end{array}
 \right] \notag\\
 \vphantom{\int^{A^A}}
 & =
  \frac{1}{u_+ - u_-}
  \left[
\begin{array}{cc}
 u_+ e^{\si_+ t} - u_ -e^{\si_- t} &  u_+ u_- (e^{\si_- t} -   e^{\si_+ t})\\
 e^{\si_+ t} -  e^{\si_- t}   &   u_+ e^{\si_- t} - u_- e^{\si_+ t}
 \end{array}
 \right], \qquad t\ge 0.
 \label{M}
\end{align}

\section{The eventual composition of the population}
\label{Sec_4}

As $\si_+ > \si_-,$ it is clear from~\eqref{M} that  $\si_+$ is the so-called {\em
Malthusian parameter\/} of the branching process that determines the long-term
behaviour of the process means. As we said in Section~\ref{Sec_1}, the most interesting
for us is the case of {\em supercritical processes\/} for which $\si_+
>0$, implying unbounded exponential growth of the population (unless it becomes extinct
at a pretty early stage). Otherwise, the process $\{\bZ (t)\}$ would be doomed to die
out very soon, so that no epidemic would arise.

It is clear that, in the supercritical case, the ratio of the time $t$ expected number
of type 2 particles to that of type~1 particles   will be given by
\[
 R_1 (t)
 \equiv \frac{M_{12} (t)}{M_{11} (t)}
 =  \frac{u_+ u_- (e^{\si_- t} -   e^{\si_+ t})}{u_+ e^{\si_+ t} - u_ -e^{\si_- t}}
 = |u_-|\biggl[1- \Bigl(1- \frac{u_-}{ u_+} + o(1)\Bigr)e^{-\D t} \biggr], \quad t\to\infty,
\]
if the process starts with a single type~1 particle, and by
\[
 R_2 (t)
 \equiv \frac{M_{22} (t)}{M_{21} (t)}
 =  \frac{ u_+ e^{\si_- t} - u_- e^{\si_+ t}}{e^{\si_+ t} -  e^{\si_- t} }
= |u_-|\biggl[1+\Bigl(1-  \frac{u_+ }{u_-} + o(1)\Bigr)e^{-\D t}\biggr], \quad
t\to\infty,
\] 
provided that it started with a type~2 particle (recall that $u_- <0$ and
$\si_+-\si_-=\D>0$).

Therefore, using~\eqref{Means}, we see that, regardless of the initial state $(z_1,
z_2)$, the eventual ratio of the mean number of $T_2$-infected hosts to that for
$T_1$-infected hosts will be one and the same quantity
\[
 R
 \equiv \lim_{t\to\infty} \frac{\exn \bigl( Z_2 (t) |\, \bZ (0) = (z_1, z_2)\bigr)} %
{\exn \bigl( Z_1 (t) |\, \bZ (0) = (z_1, z_2)\bigr)} = |u_-|
 \equiv \frac{\g_2 - \g_1 +\D}{2 \de_2},
\]
which is a well-known fact from the theory of multitype branching processes (it
follows, for instance, from Theorem~1 on p.\,185 in~\cite{AtNe}, see also p.\,203, {\em
ibid.}). The convergence rate is clearly exponential: the remainder term decays as
$e^{-\D t}$. Note also the obvious facts that $R_1(0)=0$, $R_2(0)=\infty$ and that
$R_1(t)$ ($R_2(t)$) is an increasing (decreasing) function of~$t$ (so that always
$R_1(t) < R_2(t)$).

Thus the single value $R= R(\alpha_1, \alpha_2, \beta_1, \beta_2, \mu_1, \mu_2,
\lambda)$ completely specifies the eventual balance of the mean numbers of individuals
of different types in our supercritical process (whatever the initial values). This
reflects a much deeper result on the long-term behaviour of  $\{\bZ (t), \, t\ge 0\}$
--- namely, the fact that, with probability one, the scaled vector $ e^{-\si_+ t}\bZ
(t)$ will converge, as $t\to\infty,$ to a non-trivial random vector whose distribution
is concentrated on the ray $\{r\bv,\, r\ge 0\}$ collinear to the (positive) left
eigenvector $\bv$ of~$\bA$ corresponding to the Perron-Frobenius eigenvalue~$\si_+$
(see e.g.\ Theorem~2 on p.\,206 in~\cite{AtNe} and references therein). This implies
that convergence to $R$ holds not only for the ratio of the means, but for the random
variables $Z_2 (t)/Z_1 (t)$ as well: if we denote by $A$  the event $\{  \bZ (t)\neq
\bnull$\ for all $t>0\}$, then
\[
\lim_{t\to\infty}\frac{Z_2 (t)}{Z_1 (t)} = R\quad \mbox{on}\quad A
\]
(up to an event of probability zero). In words, this means that either the branching
process becomes extinct in finite time or the sizes of the subpopulations of individuals
of the two types grow unboundedly in such a way that their ratio tends to~$R$.

Observe also that the above shows how fast the composition of the population
will change if the encounter rate $\la$ switches to another value. Suppose that
the initial value is $\la'$. As we saw, after some (exponentially short) time,
the balance of types in the process will establish around the value
$R(\lambda')\equiv R(\alpha_1, \alpha_2, \beta_1, \beta_2, \mu_1, \mu_2,
\lambda')$. Now if the value of $\la$ quickly changes to $\lambda''$, then the
population will re-establish balance at a new level $R(\lambda'')$ --- again
exponentially fast, with the rate characterized by the new value of $\D$
(provided, of course, that the process will still be supercritical, i.e.\
$\si_+>0$ for the new value of $\la$).

\begin{figure}[ht] 
%
%
%
%
\centering
\begin{tabular}{cc}
   \includegraphics[width=6cm]{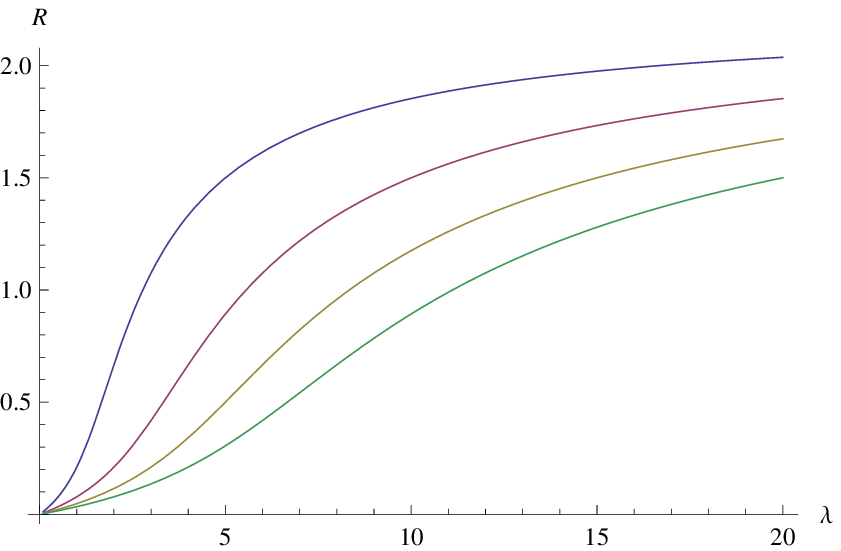} &
   \includegraphics[width=6cm]{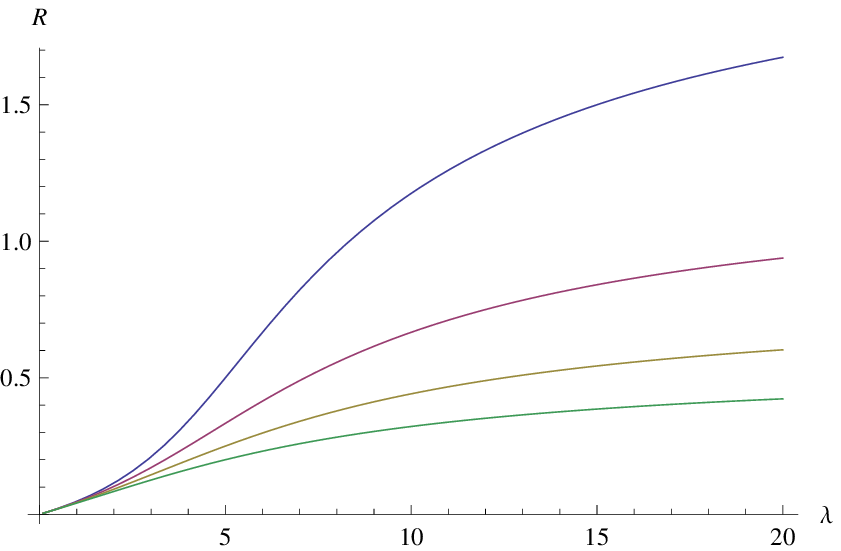}
\end{tabular}
\caption{\small The plots of $R$ as a function of $\lambda\in (0,20)$. For the fixed common values $\alpha_1=0.5,$ $\mu_1=0.2,$ $\beta_1=0.3$ and $\beta_2=0.6$,
the left pane displays the plots of $R$ for four different lethality levels $\alpha_2 =1, 1.5, 2$ and 2.5 for fixed $\mu_2=0.2$ (the lower the value of
$\alpha_2,$ the higher the curve), whereas the right one shows the plots for different mutation rates $\mu_2 =0.2, 0.3, 0.4$ and $0.6$
(the higher the mutation rate, the lower the curve), for fixed $\alpha_2=2$.
}\label{fig5}
\end{figure}

As we discussed earlier, the increase in the encounter rate $\lambda$ ought to
be beneficial for the parasite type with higher lethality as it gains more time
to spread in the host population. Indeed, we have
\[
\frac{\partial R}{\partial \la}
 = \frac{(\a_2 - \a_1)(\g_2 - \g_1 + \D)}{2\de_2 \D}>0
\]
since $\a_2 > \a_1$ by assumption and it is obvious that $|\g_2 - \g_1| <\D$.
Fig.~\ref{fig5} displays the dependence of $R$ on  $\lambda$ varying in (0,20),
for different levels of the lethality $\alpha_2$ (left pane) and mutation rate
$\mu_2$ (on the right) of the type~2 pathogen.

This example demonstrates that our model is capable of reproducing the growth
of the frequency of more lethal parasites in a host population when the density
of the hosts increases. Observe that the threshold value  $R=1$ (after which
type~2 pathogen dominates the population) may play no critical role: our model
is a crude approximation for the initial stage of an epidemic only, so
$R/(R+1)$ will just give a proportion of the carriers of type~2 parasite in the
population  of the infected hosts at the end of that stage. Even a relatively
law value of that quantity may be fatal for the population.

It is important  to note that it is the ``splitting" of the single ``virulence"
characteristic of the parasite into two (lethality and transmissibility) that
made such a capability possible: if, say, there is no difference in lethality
($\alpha_1=\alpha_2$) then, as a simple algebraic calculation shows, the value
of $R$ does not depend on density~$\lambda$. The last observation we could have
actually made earlier, as it follows from the model construction.

\begin{figure}[ht] 
%
%
%
%
\centering
\begin{tabular}{cc}
   \includegraphics[width=5cm]{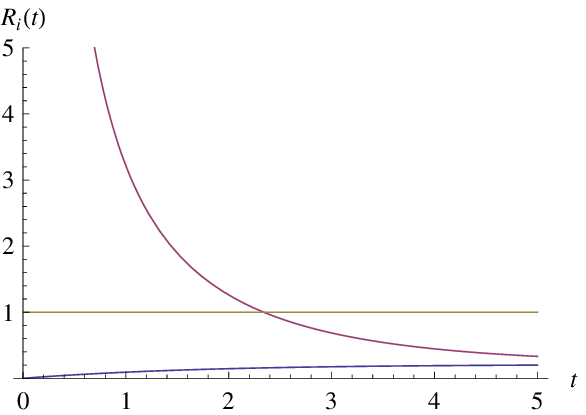} &
   \includegraphics[width=5cm]{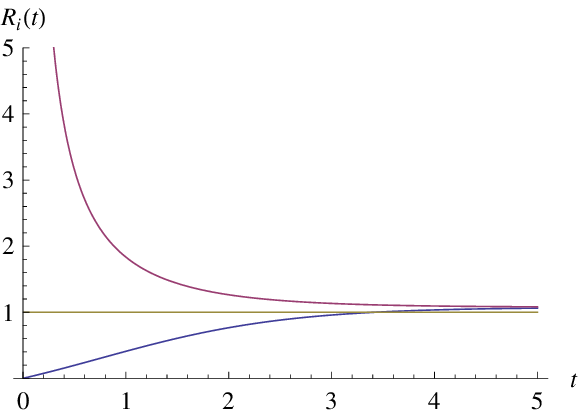}\\
   \includegraphics[width=5cm]{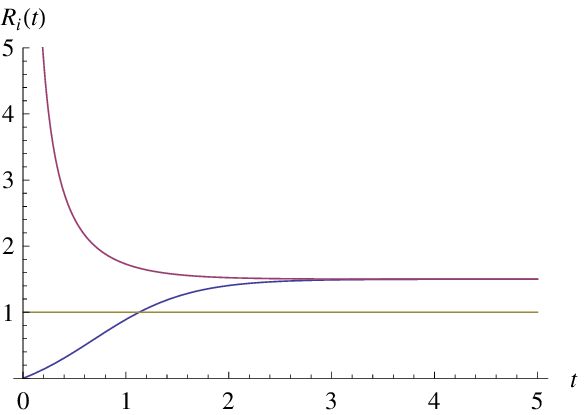} &
   \includegraphics[width=5cm]{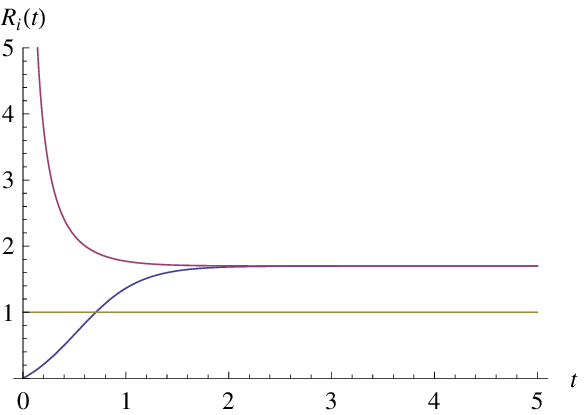}
\end{tabular}
 \caption{\small Convergence of $R_i(t)$ to $R$ as $t\to\infty$. For a common set of
parameter values, the plots display the behaviour of $R_i(t),$ $i=1,2,$ for $\la =2, 6, 10$ and $14$
(from left to right, top to bottom), with the respective values $R\approx 0.210,$ 1.076, 1.500 and 1.699. In all the cases, the process is supercritical ($\si_+>0$).}
\label{fig2}
\end{figure}

Figure~\ref{fig2} illustrates the stated exponentially fast convergence of the
ratios $R_i(t)$ to a common limit~$R$ in four situations that have common
parameter values $\alpha_1 = 0.5$, $ \alpha_2 = 1.5$, $ \mu_1 = \mu_2 = 0.2$,
$\beta_1 = 0.3$, $\beta_2 = 0.6$, but different encounter rates. The plot in
the top left corner displays the graphs of $R_1(t) < R_2(t)$ in the case when
the encounter rate $\la  = 2$ is small enough to allow type~1 parasite to
dominate ($R\approx 0.210$). On the top right plot, we see that, for $\la =6$,
there establishes a rough balance ($R\approx 1.076$), whereas on the plots in
the second raw we see type~2 parasites to gain dominance pretty fast (which is
due to higher values of $\Delta = \sigma_+-\sigma_-$), with the limiting values
$R\approx 1.500$ and $1.699$, resp.

\begin{figure}[ht] 
%
%
%
%
\centering
\begin{tabular}{cc}
   \includegraphics[width=6cm]{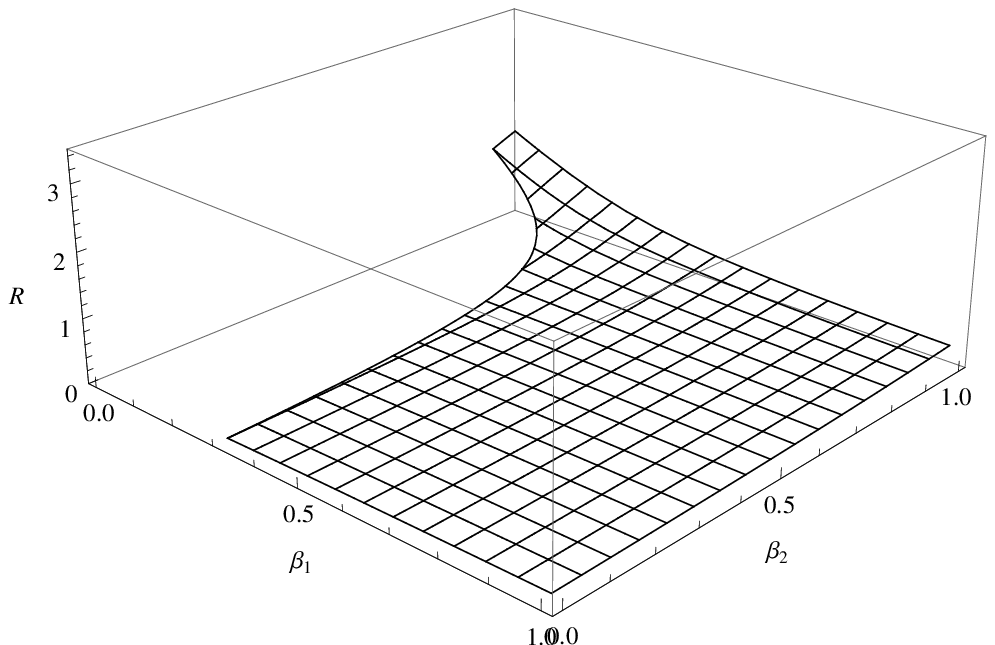} &
   \includegraphics[width=6cm]{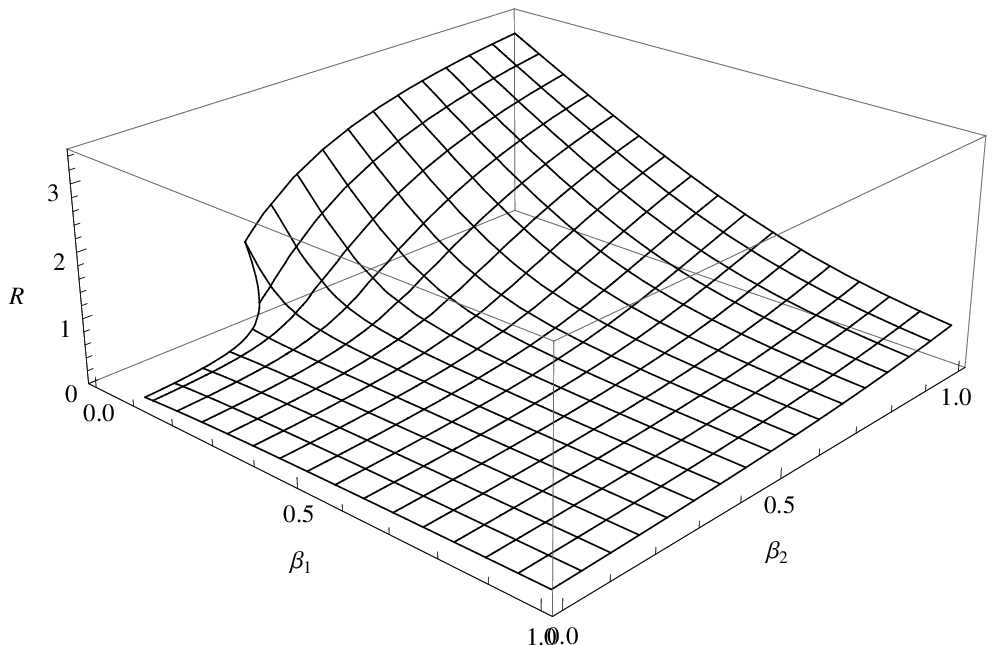}\\
   \includegraphics[width=6cm]{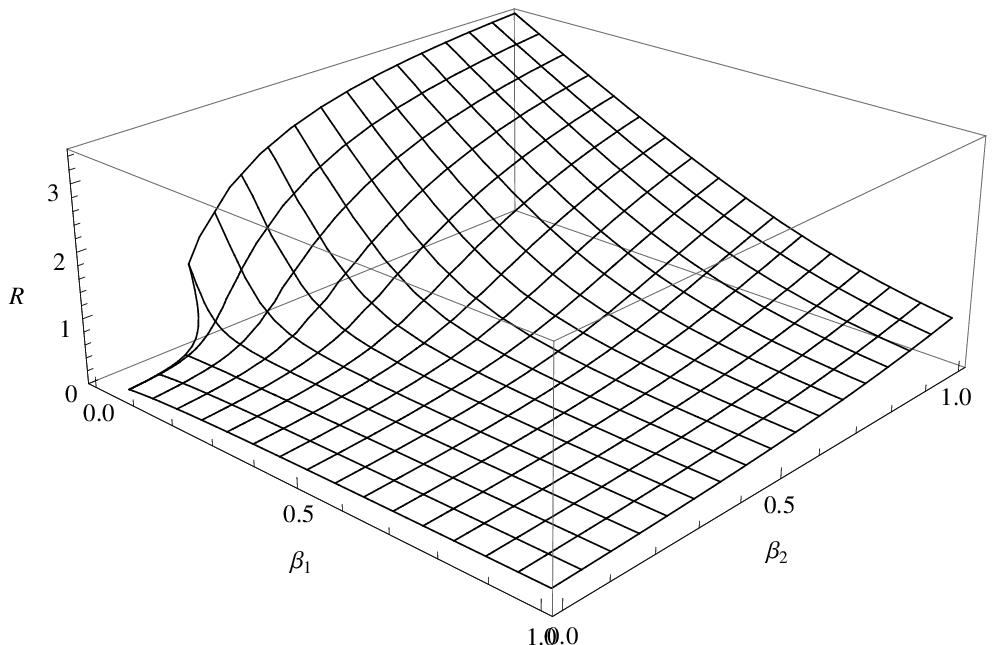} &
   \includegraphics[width=6cm]{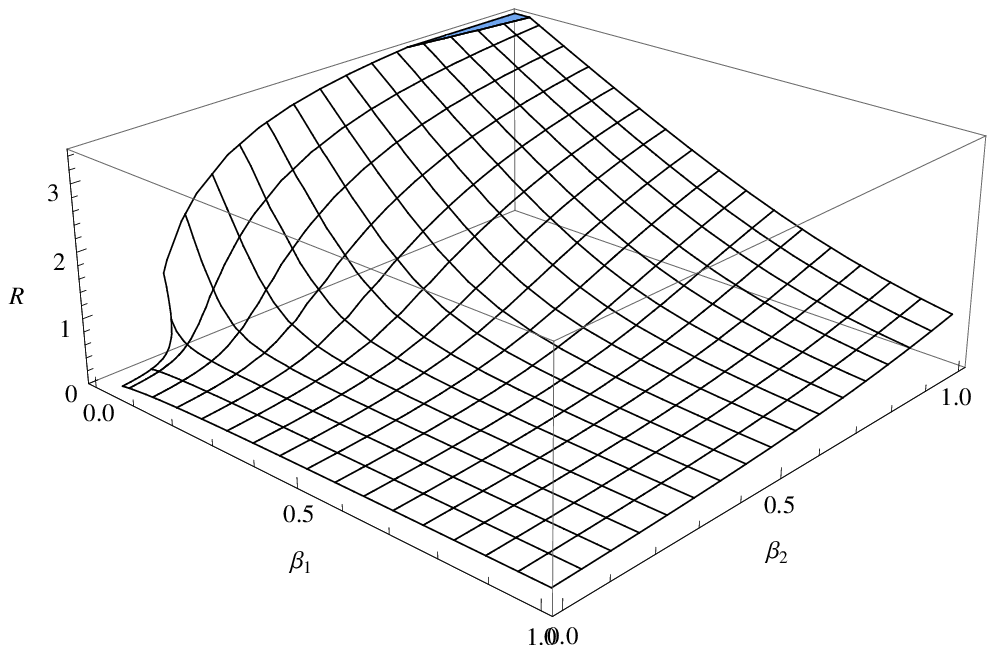}
\end{tabular}
 \caption{\small The plots of the limiting value $R$ as a function of the transmissibility probabilities $(\beta_1, \beta_2)$ in
four cases: $\la = 2, 6,10$ and 14 (from left to right, top to bottom), when all the other parameters of the model are common
($\alpha_1 = 0.5$, $ \alpha_2 = 1.5$, $ \mu_1 = \mu_2 = 0.2$). The plots are restricted to the regions where the respective processes are supercritical.
 }

\label{fig3}
\end{figure}

The character of the dependence of $R$ on the transmission probabilities is
illustrated in Fig.~\ref{fig3}. For   four different values of the encounter
rate ($\lambda=2,6,10$ and~14), the figure shows the plots of $R$ as a function
of $(\beta_1, \beta_2)$, restricted to the regions where the process is
supercritical (i.e.\ $\sigma_+>0$). The values of the other parameters are
$\alpha_1 = 0.5$, $ \alpha_2 = 1.5$, $ \mu_1 = \mu_2 = 0.2$. As one could
expect, the value of $R$ strongly depends on $(\beta_1, \beta_2)$ and is an
increasing function of~$\beta_2$ and a decreasing one of~$\beta_1.$

\begin{figure}[ht] 
%
%
%
\centering
\begin{tabular}{cc}
   \includegraphics[width=6cm]{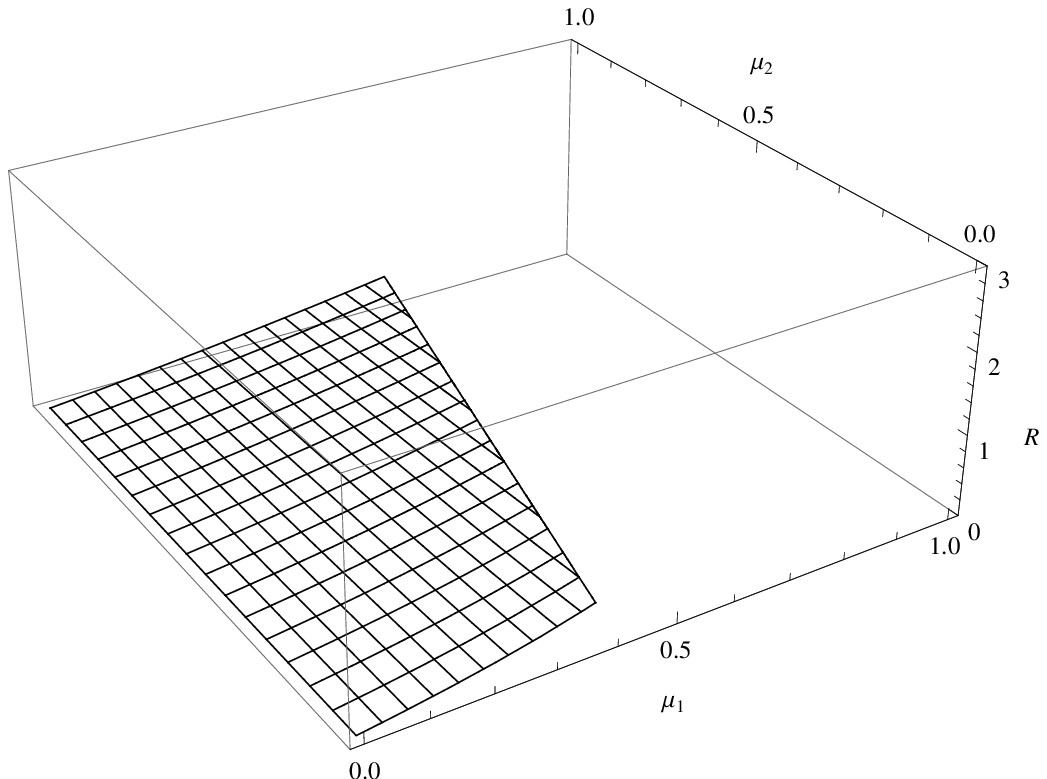} &
   \includegraphics[width=6cm]{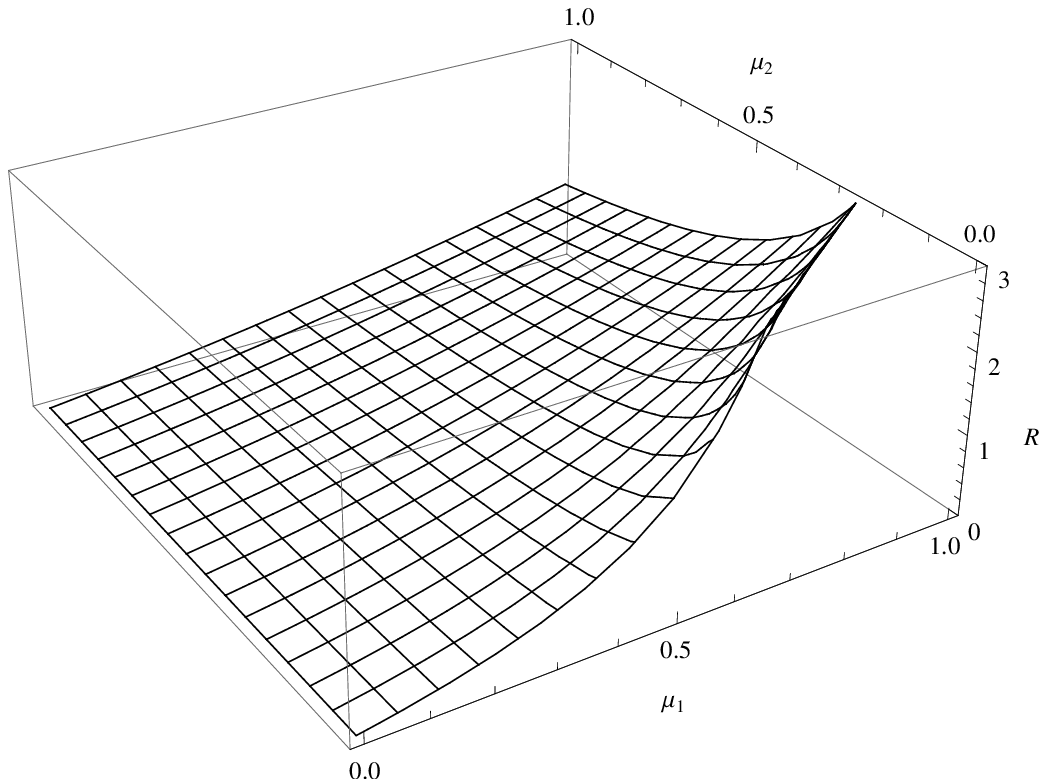}
\end{tabular}
\caption{\small The plots of $R$ as a function of the mutation probabilities for $\la=4$ (the left pane) and $\la=7$ (on the right), for fixed common values of
all the other parameters of the model ($\a_1=1,$  $\a_2=2,$  $\b_1=0.4$ and $\b_2=0.2.$ The plots are restricted to the regions where the respective processes are supercritical.
Note the unusual orientation of the $\mu$-coordinate axes
(chosen so as to have a better view of the plots).
}\label{fig3+}
\end{figure}

In all the four cases presented in Fig.~\ref{fig3} the threshold value $R=1$ is
exceeded only when the transmissibility of type~2 pathogen is greater than that
for type~1 ($\beta_2 > \beta_1$), so it may appear that that inequality is a
necessary condition for type~2 to   prevail. This, however, is not true: it
turns out that a higher mutation rate  from type~1 to type~2 can compensate for
some lack of transmissibility.  Figure~\ref{fig3+}   shows the plots of $R$ as
a function of the mutation probabilities $(\mu_1, \mu_2)\in (0,1)^2$ for
different encounter rates, all other parameters being fixed and common, with
the transmission probability for type~1 being double that for type~2 ($\b_1=
2\b_2=0.4$). On the left plot corresponding to $\la=4$, the maximum value of
$R$ barely exceeds 0.5, whereas on the right one, due to the increase in the
encounter rate to $\la=7,$ not only the  supercriticality region is much
bigger, but also the maximum value is~$R=3$. We see type~2 parasite's
domination in the region where the mutation rate $\mu_1$ (from type~1) is high
enough, while $\mu_2$ is relatively small.

\begin{figure}[ht] 
%
%
%
%
\centering
\begin{tabular}{cc}
   \includegraphics[width=6cm]{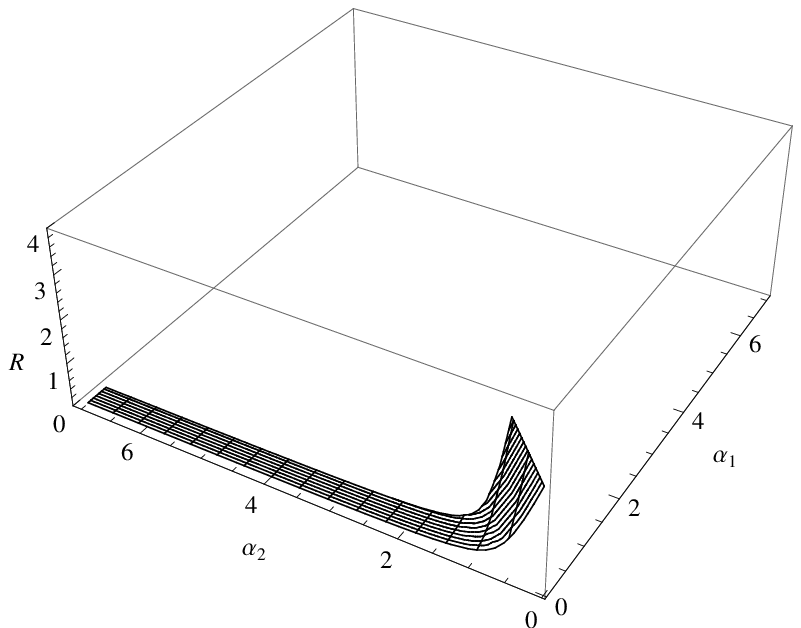} &
   \includegraphics[width=6cm]{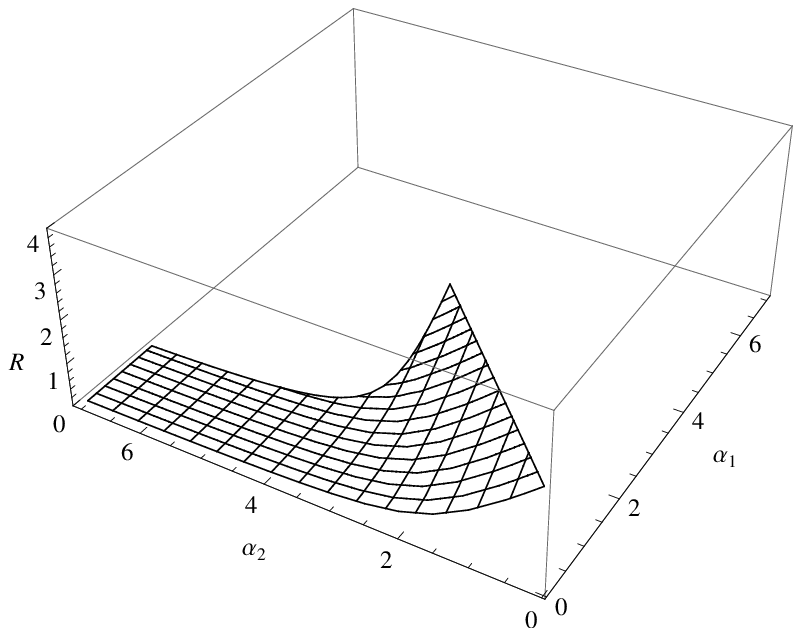}\\
   \includegraphics[width=6cm]{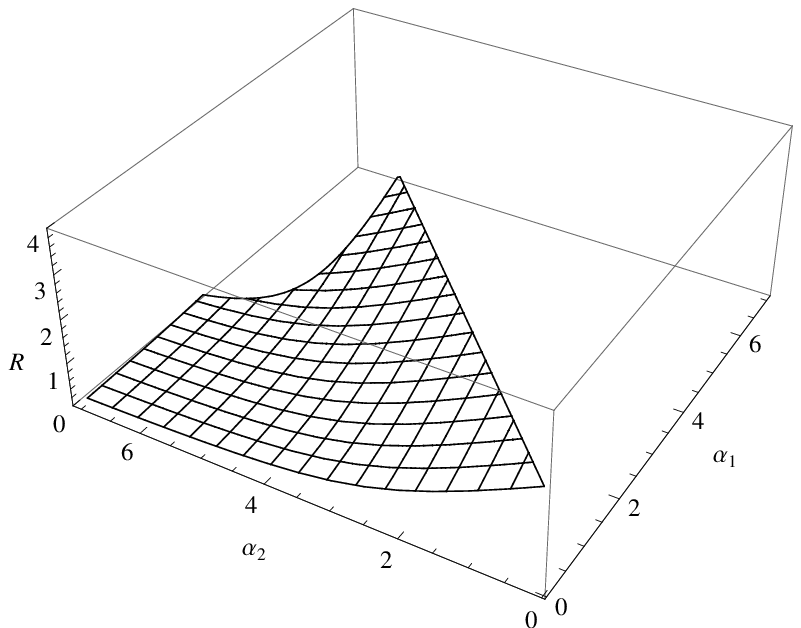} &
   \includegraphics[width=6cm]{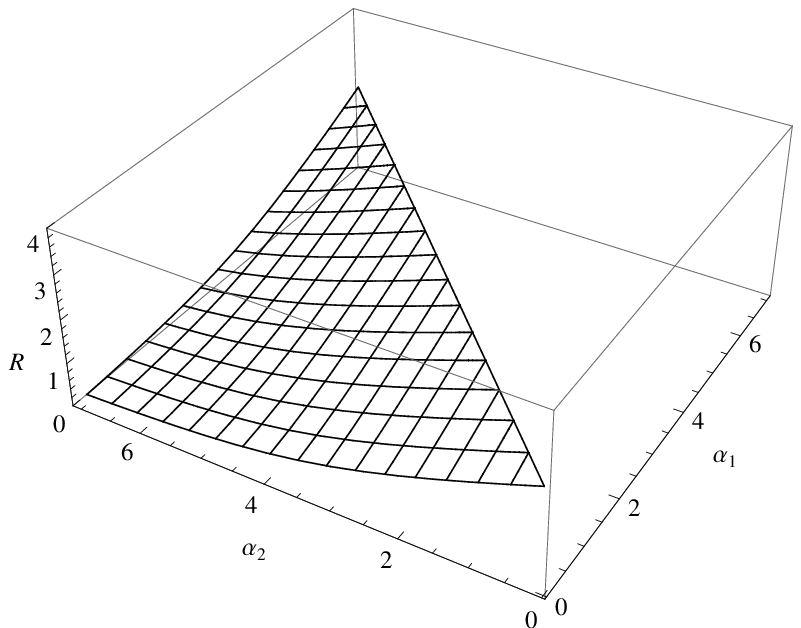}
\end{tabular}
\caption{\small
The plots of  $R$ as a function of the lethalities $(\alpha_1, \alpha_2)$, in
the same four cases as on Figs.~\ref{fig2},~\ref{fig3}: $\la = 2, 6,10$ and 14
(from left to right, top to bottom), when all the other parameters of the model
are common ($\beta_1 = 0.3$, $ \beta_2 = 0.6$, $ \mu_1 = \mu_2 = 0.2$). The
plots are restricted to the regions where the respective processes are
supercritical and $\alpha_1 <\alpha_2$ (as we assumed). Note the unusual
orientation of the $\alpha$-coordinate axes (chosen so as to have a better view
of the plots).}
 \label{fig4}
\end{figure}

Finally, we turn to the dependence of $R$ on the lethalities $\alpha_j$. As one
can easily see,
\[
\frac{\partial R}{\partial \a_1} = \frac1{2\de_2}\Bigl(1 - \frac{\g_1 - \g_2}{\D}\Bigr)<0,
 \qquad
 \frac{\partial R}{\partial \a_2} = - \frac{\partial R}{\partial \a_1}>0.
\]
This is quite natural, as the increase in a pathogene type's lethality does not
improve its chances  to prevail when all other parameters of the model remain
unchanged. The character of the dependence is illustrated in Fig.~\ref{fig4}
showing $R$ as a function of $(\alpha_1, \alpha_2)\in (0,7)^2,$ for $\la=
2,6,10$ and~14. Note how the encounter rate $\la$ influences the size of the
region where the process is supercritical (thus, for $\la_2$ it shrinks to a
narrow strip in the $(\alpha_1, \alpha_2)$-domain, corresponding to small
values of $\a_1$).

\section{On the effects of the change in enclosure size on the encounter rate}
\label{Sec_5}

As we said, our main motivation was to model the effect that the change in ``effective
density" represented by the parameter $\la$ (the rate at which infected hosts encounter
susceptibles) can have on the ratio of the number of hosts infected, say, with $T_2$
pathogen to that for $T_1$-infected hosts. But how does the value of $\la$ relate to
physically measurable parameters of the modelled situation\,---\,for instance, the
density of fish in an aquaculture tank (given the tank size and all other parameters of
the model are fixed) or the size of the tank (given the number of fish in it is fixed)?
How will $\la$ change if one, for example, ``squashes" the same host population to a
``world" whose linear dimensions are twice as small as for the original one?

The answer to this type of questions will in general depend on what one assumes about
the character of the hosts' movements (and of course, on the the pathogen transmission
mechanism\,---\,but we will not address this aspect in our simple analysis in the
section). One of the most popular models for ``wandering particles" is the famous
Brownian motion process $\{W(t), t\ge 0\}$ (see e.g.\ p.169 in~\cite{KaTa}), which can
be thought of as a continuous analog of a simple (symmetric) random walk. Recall that
the Brownian motion is defined as a continuous time process with continuous
trajectories that starts at zero at time $t=0$ and has independent Gaussian increments:
$W(t+h) - W(t) \sim N(0,h)$ for $t,h\ge 0$. One of the key properties of the process is
its {\em self-similarity\/}: for any $a>0,$ one has
\begin{equation}
\{a  W(t), t\ge 0\}\stackrel{d}{=} \{W(a^2 t), t\ge 0\} ,
 \label{SS}
\end{equation}
i.e.\ these two processes have the same distribution.

Since the total number of hosts is assumed to be large enough, encounter rates are
mostly determined by the ``local" characteristics (the density of the population and
the dimensionality of the space) of the enclosure and will have little dependence on
the ``shape of the world". Therefore, for analysis purposes, we will assume in this
section that hosts perform independent Brownian motions in a ``simple world" $S$ in
one, two or three dimensions (starting at some ``individual" initial points) and show
how the encounter rates change when one ``contracts" the original world without
changing its shape, i.e.\ switches from $S$ to the set $\ep S=\{\ep \bx:\, \bx\in S\}$.
In this context, the dimensionality can actually be though of as a crude description of
the the shape of the world.

Suppose that there are $N$ susceptibles in the population and that the movements of all
the hosts are independent of each other. The value of the parameter $\la=\la_N$ gives
the average number of encounters of a given infected host with susceptible ones per
time unit (we assume that $n$ is large enough so that the ``conversion" of some
susceptibles into infected hosts during our modelling ``time horizon" does not
significantly change the number $N$ and hence the encounter rate $\lambda$). It is
clear that $\la_N = N\lambda_1$ and, moreover, that $\lambda_1 = 1/t^*$, where $ t^* $
is the mean time to encounter of our infected host with a given susceptible host. Thus
the answer to the question on how the encounter rate will increase if the host density
living in a "fixed world" increases is simple: it is just proportional to the number of
hosts in the world. However, to understand the effect of the {\em world size change\/}
(when we switch from~$S$ to~$\ep S$) on the encounter rate $\la = \la (\ep)$, we will
have to analyze that effect on the mean time~$t^* = t^*(\ep)$.

\medskip\noindent
{\bf The one-dimensional case.} First we assume that $S = [0,1]\subset \mathbb{R}$. Our
hosts move according to independent Brownian motions, reflecting from the boundaries~0
and~1 of the set~$S$. This can be formalised by introducing the function
\[
\varphi (x) := \left\{
\begin{array}{ll}
x - \lfloor x\rfloor &\mbox{if $\lfloor x\rfloor$ is even},\\
1 -x + \lfloor x\rfloor &\mbox{if $\lfloor x\rfloor$ is odd},\\
\end{array}
\right.  \qquad \lfloor x\rfloor := \max\{k\in \mathbb{Z}:\, k\le x\},
\]
and letting the location of our infected host to be given by $H_0=H_0(t) = \varphi (H_0
(0) + W_0 (t))$ and that of a given susceptible one by $H_1=H_1(t) = \varphi (H_1 (0) +
W_1 (t))$, where $H_i (0)\in S$ are some fixed initial positions and $W_i$ are
independent standard Brownian motions, $i=0,1$. We say that the hosts have an encounter
at time $t$ if $H_0(t) = H_1(t)$.

Now consider the space $\ep S$, where our hosts move according to $\{\ep H_i (t),\,
t\ge 0\}$ and denote by $t^* (\ep)$ the mean time to encounter of the hosts in this new
world, $\ep >0.$ It is obvious from the self-similarity property~\eqref{SS} that $t^*
(\ep ) = \ep^2 t^* (1) $ and therefore
\[
\la (\ep ) = \ep^{-2} \la (1).
\]
Thus, if our hosts are confined to a one-dimensional world where they wander at random,
the encounter rate displays inverse quadratic dependence on the world size: say, halving
the ``living space" will increase the encounter rate fourfold.

\medskip\noindent
{\bf The two-dimensional case.} To avoid dealing with any boundaries, we assume that
our hosts live on the two-dimensional sphere
$$
\mathbb{S}^2 = \{ (x,y,z) \in \mathbb{R}^3:\, x^2 + y^2 +z^2 =1\}
$$
and that our hosts wander on it at random according to independent spherical Brownian
motions $\{H_i(t),\, t\ge 0\}$  (see e.g.\ \cite{KaTa}, Chapter~15, Section~13I),
starting at some fixed distinct points $H_i(0),$ $i=0,1$.

In this case, $\Pr (H_0 (t)\neq H_1 (t),\, t\ge 0)=1$, so we need to modify our
definition of encounter. Fix a small enough $\de >0$ and define the encounter time of
the hosts $H_i$ as $\inf\{t>0: r(H_0(t), H_1 (t))=\de\}$, where $r(\cdot, \cdot)$ is the
geodesic distance on~$\mathbb{S}^2$. Denote by~$t^*_\de = t^*_\de (1)$ the mean value of
this time (suppressing the dependence of the initial locations $H_i(0)$; to simplify the
exposition, we deliberately make it somewhat sketchy).

Next we consider the ``contracted world" $\ep \mathbb{S}^2$, where our hosts wander
according to $\{\ep H_i(t),\, t\ge 0\}$, but the definition of encounter remains
unchanged (the hosts should find themselves within distance~$\de$ of each other); the
mean time to encounter for this case is denoted by $t^*_\de (\ep)$. Again using
self-similarity, we can easily conclude that
\begin{equation}
t^*_\de (\ep ) = \ep^2 t^*_{\de/\ep} (1) .
 \label{tt}
\end{equation}
However, we want to relate $t^*_\de (\ep )$ to $ t^*_{\de} (1),$ so it remains to
clarify the relationship between $t^*_{\de/\ep} (1) $ and $t^*_{\de} (1)$.

It is obvious that $t^*_{\eta} = t^*_{\eta} (1) $ is a decreasing function of $\eta>0$.
As we are interested in situations where $\de/\ep$ is small (despite the small size of
the ``world", encounters are still relatively rare), it suffices to find the asymptotic
behaviour of $t^*_{\eta}  $ as $\eta \to 0$. To do that, we first observe that analyzing
the dynamics of the position of $H_0 (t)$ relative to $H_1(t)$ shows that finding the
mean time when the two points are first within distance~$\eta$ of each other is
equivalent to finding the mean time a Brownian particle $H^*(t)$ (with an initial
position at a distance~$r(H_0(0), H_1(0))$ from the ``North Pole" $P=(0,0,1)\in
\mathbb{S}^2$ and local diffusion coefficient $\sqrt{2}$ times that for the original
spherical Brownian processes) will need to get within distance $\eta$ of $P$.

Denoting by $V(t)$ the projection of $H^* (t)$ on the $z$-axis, one can easily see from
It\^o's formula that $\{V(t), \, t\ge 0\}$ is a diffusion process with the state space
$[-1,1]$ governed by the stochastic differential equation
\[
dV(t) = \mu (V(t)) dt + \sigma (t) dW_0(t),
\]
where $W_0$ is a standard univariate Brownian motion process and the drift and diffusion
coefficients are given by
\[
\mu (z) = - \frac{z}2, \qquad \sigma^2(z) = \frac{1-z^2}2,
\]
respectively (see e.g.~p.194 in~\cite{Ma}; for convenience we assumed that $H^*(t)$
follows a standard spherical Brownian motion on $\mathbb{S}^2$ which will have no
adverse implications for the validity of our analysis). The geodesic distance from
$H^*$ to $P$ is equal to $\eta$ iff its projection on the $z$-axis equals $r:=\cos
\eta,$ so that we need to find
\[
v(z) = \exn (\tau |\, V(0) =z), \qquad\mbox{where}\quad \tau= \inf\{t>0: \, V(t) =r\}.
\]
This can easily be done using the standard technique of the method of differential
equations (see e.g.\  Problem~B on p.192 in~\cite{KaTa}): the function $v(z)$ is the
bounded solution to the equation
\begin{equation}
1 = -\mu (z) v'(z) - \frac12 \sigma^2 (z) v'' (z)
 \equiv  \frac12 zv'(z) + \frac14 (z^2 -1) v'' (z), \qquad z\in (-1,r),
\label{DE}
\end{equation}
with the boundary condition $v(r)= 0$. Setting $u(z): =  (z^2-1)/4 ,$ we notice that
$u'(z) =z/2$ and so \eqref{DE} is equivalent to
\[
1 = u'v' + uv'' = (uv')',
\]
which means that $u(z) v'(z) = z+c_1,$ and so
\[
v'(z) = \frac{z+ c_1}{u(z)} = \frac{4(z+c_1)}{z^2-1}.
\]
Therefore the general solution to \eqref{DE} is given by
\[
v(z) = 2\bigl[(1+c_1) \ln (1-z) + (1-c_1) \ln (1+z)\bigr]+ c_2,
\]
which is bounded on $(-1,r)$ iff $c_1=1.$ Now using the boundary condition at $z=r$ to
find $c_2$ leads to
\[
v(z) = 4\ln\frac{1-z}{1-r}, \qquad   z\in [-1,r].
\]
To find the asymptotic behaviour of $v(z)$ as $\eta\to 0$, we use $\cos \eta= 1 -
\eta^2/2 (1+o(1))$ to obtain that, for a fixed initial value~$z$,
\[
v(z) = 8 |\ln \eta| + O(1).
\]

This means that, for fixed initial positions of $H_0$ and $H_1$, we have $t^*_\eta =
(c+o(1)) |\ln \eta|$ as $\eta\to 0$. Together with \eqref{tt} this yields, for small
$\de/\ep,$
\[
t^*_\de (\ep )\approx \ep^2 \frac{|\ln \de|- |\ln \ep|}{|\ln \de|}  t^*_{\de} (1).
\]
That is, the encounter rate behaves as
\[
\la (\ep) \approx \ep^{-2} \frac{|\ln \de|}{|\ln \de|- |\ln \ep|} \la (1).
\]
In the two-dimensional case, it might be more natural to relate the encounter rate not
to the linear dimensions of the enclosure, but rather to its area (proportional
to~$\ep^2$). The above formula shows that, in this case, the encounter rate in a
``shrinking" world still grows somewhat faster than the inverse proportional to its
area, the latter giving the density of hosts.

\medskip\noindent
{\bf The three-dimensional case.} We still have \eqref{tt}, and an analysis similar to
the one carried out in the two-dimensional case shows that now $t^*_\eta = (c+o(1))
\eta^{-1}$ as $\eta\to 0$ (cf.\ p.195 in~\cite{Ma}), so that $ t^*_\de (\ep ) = \ep^3
t^*_{\de} (1)(1+o(1)) $. That is,
\[
\la (\ep) =   \ep^{-3}   \la (1) (1+o(1)).
\]
We see that, in the three-dimensional case (assuming that the hosts wander according to
three-dimensional Brownian motions), the encounter rate in a ``shrinking" world is
inverse proportional to its volume. That is, in this case the rate is proportional to
the density of hosts.

\medskip To summarise the above analysis, we observe that the encounter rate $\la=\la (\ep)$
grows rather fast when the linear dimensions (specified by the parameter~$\ep$)
of the hosts' "world" diminish. In the three-dimensional case, $\la$ is
inversely proportional to the volume per host, in the two dimensional case it
grows slightly faster than the reciprocal of the area per host, while in the
one-dimensional case the growth rate of $\la$ is inversely proportional to the
square of the size of a host's share of the enclosure. This indicates that not
only effective density per se, but also the shape of the enclosure  can be an
important factor leading to an epidemic. Thus the nature of the enclosures in
which animals are kept can be an important factor in determining the progress
and nature of an epidemic.

\section{A multistage modification of the model}
\label{Sec_6}

In this section we will consider an aggregate model for situations in which there are
several populations of hosts that exist in originally isolated enclosures.

There may be many pens with animals at a farm or many fish tanks at an aquaculture
facility. Initially, one of the enclosures becomes infected with a single type
pathogen. This can give rise to a ``local epidemic" in the infected enclosure, which
can be modelled using our processes from Section~\ref{Sec_2} (assuming that we have
supercriticality: $\si_+>0$). The original pathogen may also mutate to become more or
less lethal. We will initially assume that it may mutate to a more lethal type~2
pathogen ($\alpha_2 =r \alpha $ for some $r>1$, $\alpha_1 =\alpha$). That new pathogen
type can also have different transmissibility and mutation rate, but, to make our model
as simple as possible, we will assume for the time being that it differs from type~1
pathogen in lethality only, all other parameters being common. Denoting them simply
by~$\beta$ and $\mu$, we see that the Malthusian parameter for that process is given by
\begin{align}\label{siGaa}
\si_+(\alpha_1, \alpha_2 ) = \frac12\left( 2 (1-\mu)\beta \la- \alpha_1 - \alpha_2
+\sqrt{(\alpha_1 - \alpha_2)^2 + 4 \mu^2\beta^2\la^2}\right).
\end{align}

After the epidemic has gone through the initial stage (and so the ratio of the numbers
of  hosts infected with different pathogen type can be assumed equal to $\rho_0$, where
$\rho_k \equiv R(r^k \alpha , r^{k+1} \alpha , \beta , \beta , \mu , \mu , \lambda)$),
the infection is transmitted to the next enclosure (say, by a worker in a farm
situation). The transmitted pathogen is chosen at random, so that the probability of
transmitting the one with lethality $\alpha$ (denote this event by~$A$) is
$1/(1+\rho_0)$, while the one with lethality $r \alpha$ is transmitted with probability
$\rho_0/(1+\rho_0)$. This, in turn, may lead to a local epidemic in the new enclosure:
we again assume the possibility of mutation to a more lethal pathogen (so that now we
will have $\alpha_1 = \alpha,$ $\alpha_2 = r\alpha$ if the event~$A$ occurred, and
$\alpha_1 = r\alpha,$ $\alpha_2 = r^2\alpha$ otherwise), and to have an epidemic we
again need $ \si_+ = \si_+(\alpha_1, \alpha_2 ) >0 $ (now for the new set of parameters
$\alpha_1,$ $ \alpha_2$). Once the epidemic has established itself in the second
enclosure (and the balance of pathogen types has stabilized around the respective
$R$-value), the next step is the transmission of the disease (by means of a random
mechanism of the same type as in the first instance) to the next enclosure, and so on.

Scenarios of this type have been encountered often where once a disease is recognized
in a herd, animals in the infected enclosure are removed or killed, but the disease is
subsequently found in other herds, for example the spread of Foot and mouth disease
among herds in Taiwan, which was related to herd size and the number of herds in a
province~\cite{Ge}. Of course, biosecurity measures are intended to prevent such
transmission between enclosures, but often the need for diligence is learned after the
event.



It is easily seen from~\eqref{siGaa} that
\begin{equation}
 \label{deriv}
\si_+(0,0)>0,\qquad
 \frac\partial{\partial \alpha} \si_+( \alpha , r \alpha )<0,
 \qquad\lim_{\alpha \to\infty}  \si_+( \alpha , r \alpha ) =-\infty,
\end{equation}
and also that $\dfrac\partial{\partial \alpha} R ( \alpha , r \alpha )<0$. Thus, if the
lethality of the pathogen will keep increasing, the Malthusian parameter of the model
will eventually drop below zero, and then the epidemic will collapse. More
specifically, setting
\[
k^* \equiv \inf\{k\ge 0:\, \si_+(r^k \alpha , r^{k+1}\alpha )<0\},
\]
we see that
\begin{equation}
 \label{Maltneg}
 \si_+(r^k \alpha , r^{k+1}\alpha )<0 \quad\mbox{ for all $k\ge k^*$.}
\end{equation}


It is clear that the transition of the disease from enclosure  to enclosure
according to the above scheme is described by a discrete time Markov chain
$\{X_n\}$, the ``time" $n$ having the meaning of the number of steps (i.e.\
enclosures infected), $X_n$ representing the level of lethality of the
pathogens in the $n$th infected enclosure: we set $X_n=k$ if
$(\alpha_1,\alpha_2) =(r^k \alpha , r^{k+1} \alpha)$ in the enclosure. Thus the
state space of the Markov chain is $\{0,1,2,\dots\}$ and the only nonzero
entries in the transition matrix $\bP = [p_{j,k}]_{j,k\ge 0}$ of the chain are
\begin{align*}
 p_{k,k} &\equiv \Pr (X_{n+1}= k |\, X_{n+1}= k) = \frac1{1+\rho_k}, \\
 p_{k,  k+1} &\equiv \Pr (X_{n+1}= k+1 |\, X_{n+1}= k) = \frac{\rho_k}{1+\rho_k},
\end{align*}
$k=0,1,2,\dots$

Further, in view of~\eqref{Maltneg}, one can assume that once the Markov chain
$\{X_n\}$ has reached the state $k^*,$ the epidemic becomes unsustainable, and
hence there will be no further transmission of the disease to other enclosures.
So we can truncate our state space to $\{0,1,2,\dots, k^*\}$, which results in
a finite decomposable Markov chain with a single absorbing state~$k^*$.
Whatever the current state of the chain, at the next step it can either stay at
it or move to the right, the transition probabilities forming the matrix
\[
[q_{jk} ]= \bQ = \left[
\begin{array}{cc}
\bT & \br^\top \\
\bnull & 1
\end{array}
\right],
\]
where $\bT$ is the $k^*\times k^*$ substochastic matrix formed by the first
$k^*$ rows and $k^*$ columns of $\bP$, $\br = (0,\dots, 0, p_{k^*-1, k^*} ) \in
\mathbb{R}_+^{k^*}$ and $^\top$ denotes transposition.

The (random) number of steps $T$ the chain will need to reach the absorbing
state is nothing else but the total number of enclosures that will be affected
by the epidemic prior to its collapse. Using our model, we can easily find the
distribution of~$T$.

Indeed, using the standard approach to solving such problems (see e.g.\ p.80
in~\cite{Ki}), we note that as the state $k^*$ is absorbing, we have $\Pr (T\le
n |\, X_0=0)= q^{(n)}_{0,k^*}$, where $q_{jk}^{(n)}$ are the $n$-step
transition probabilities:
\[
[q_{jk}^{(n)}]=\bQ^n   = \left[
\begin{array}{cc}
\bT^n & \br^\top_n \\
\bnull & 1
\end{array}
\right], \qquad \br_n^\top = (\bI + \bT + \cdots + \bT^{n-1})\br^\top,
\]
so that for the probability mass vector function $\bff (n) = \{f_j (n) , \, j=0,\dots,
k^*-1\}$, $f_j (n)= \Pr (T=n|\, X_0=j)$, one obtains
\[
\bff (n) = \br_n - \br_{n-1} = \br (\bT^{n-1})^\top , \qquad n\ge 1.
\]
Of course, we are only interested in the first entry of the vector $\bff (n)$.

To compute the mean and higher moments of $T$ one can use the generating function
\[
\bff^* (z) \equiv \sum_{n=1}^\infty z^n \bff (n)
 = \br z \bigl(\bI - z \bT^\top\bigr)^{-1}, \qquad |z| \le 1.
\]
In particular, since $\dfrac{d}{dz} \bff^* (z)= \br \bigl( \bI -  z \bT^\top\bigr)^{-1}
+ \br z \bT^\top  \bigl(\bI -   z\bT^\top\bigr)^{-2}$, we find that
\[
\exn (T |\, X_0=j) = \frac{d}{dz} \bff^* (z)\bigg|_{z=1}
 = \br \bigl(\bI -   \bT^\top\bigr)^{-2}
 = \bun \bigl(\bI -   \bT^\top\bigr)^{-1},
\]
where the last equality follows from the obvious observation that $\br   \bigl(\bI -
\bT^\top\bigr)^{-1}= \bff^* (1) = \bun\equiv (1,\dots, 1)\in \mathbb{R}_+^{k^*}$.

A possible objection to the above simple aggregate model is that pathogens will
not always mutate to become more lethal. The model can be further generalized
by allowing, within each enclosure, mutations of our pathogen not only in the
direction of higher lethality, but also in the opposite direction. So we will
first have to generalize our basic model from Section~\ref{Sec_2} to a
three-type branching process, assuming that, if an enclosure is infected with a
pathogen with lethality~$\alpha$, then the pathogen can mutate to ones with
lethalities $r^{-1}\alpha$ and $r\alpha$, respectively, where, as before, $r>1$
is a fixed number (all other parameters being assumed equal for the different
types of pathogens). Mathematically, analyzing such processes is essentially
equivalent to what we did in Sections~\ref{Sec_2}--\ref{Sec_4}, although all
the closed-from expressions will become much more complicated, and so we will
not give much technical detail for brevity's sake. The main assertions
concerning the asymptotic behaviour of the branching process will remain true:
there will exist a limiting balance of types in the supercritical case (denote
the shares of the different pathogens by $\pi_j= \pi_j (\alpha)$, $j=1,2,3,$
$\sum_j \pi_j=1$), which can be found from the generator $\bA$ of the semigroup
of the mean matrices, and the almost sure convergence of the process scaled by
$e^{-\sigma_+ t}$ (as before, $\sigma_+$ denotes the Perron-Frobenius root of
$\bA$) to a limiting random vector will hold.

In the multi-stage model, we start with initial infection of one of the
enclosures with a pathogen with lethality~$\alpha$. That leads to an epidemic
(provided, of course, that $\sigma_+>0$) in which pathogens of three types will
be present, with lethalities given by the vector  $(\alpha_1,
\alpha_2,\alpha_3)=(r^{-1}\alpha ,  \alpha ,r \alpha )$. The next enclosure to
be infected will receive a pathogen chosen at random from those present in the
first infected enclosure, and so it will have lethality $\alpha_j$ with
probability $\pi_j$, $j=1,2,3$, and so on.

Observe that the triplets of lethalities $(\alpha_1, \alpha_2,\alpha_3)$
characterizing the pathogens present in a given enclosure in our system will
all be of the form $(x , r x ,r^2 x )$ for some $x>0$, i.e.\ lying on a common
ray $L$ with the direction vector $(1 , r ,r^2)$. Therefore we will again have
a basically ``univariate" Markov chain $\{X_n\}$ showing what pathogens can be
present in different enclosures, $X_n=k$ meaning that the $n$th affected
enclosure was initially infected with the pathogen of lethality $r^k\alpha$,
$k\in\mathbb{Z}$ (and in this case there can also be pathogens with lethalities
$r^{k-1}\alpha$ and $r^{k+1}\alpha$ in that enclosure), assuming that $X_0=0$
(as $\alpha$ is the lethality of the pathogen that was initially introduced
into our system). The original state space for the process will
be~$\mathbb{Z}\equiv \{\dots, -1,0,1,\dots\}$, which is infinite in both
directions. At each step, the value of the chain can remain unchanged (the
interpretation being that the pathogen transmitted to the next enclosure had
the same lethality as the one with which the epidemic started in the current
one) or can either decrease or increase by one (that is, the transmitted
pathogen would have lethality values equal to $r^{-1}$ or $r$ times the current
one, respectively).

Further, , as before, one can show that $\frac\partial{\partial x} \si_+( x, r
x ,r^2x)<0,$ so that, moving along the ray~$L$ in the ``positive direction", we
will eventually enter the subcriticality region for the branching process,
where the basic reproductive number will be less than one. Therefore, at this
stage the epidemic will collapse, and hence we again can ``truncate" the state
space for $\{X_n\}$\,---\,now to  $\{\dots, -1,0,1,\dots, k^*\}$, where $k^*$
is an absorbing state that has the same meaning as above.

The variety of behaviours that such a chain can display will be somewhat wider
than for our first multi-stage model. The dynamics of the chain are determined
by the behaviour of the mean step values
\[
\exn \bigl( X_{n+1} - X_n \big|\, X_n =k ) = \pi_3 (r^{k}\alpha) - \pi_1 (r^{k}\alpha),
\qquad k<k^*.
\]
In particular, the parameters of the model can be such that the above
quantities will be negative. Then absorbtion at $k^*$ occurs with probability
less than one, while on the complement event the chain will drift away in the
negative direction, which corresponds to the disease ``fading", when the
pathogen's lethality vanishes, and so on.


\section{Discussion}
\label{Sec_7}

We will conclude with a few remarks concerning possible biological
interpretation of our results.

The mathematical models we presented show that, at the beginning stages of any
epidemics that arise in situations where animals (or humans) live in enclosures, the
{\em density of hosts\/} is an important factor in determining whether more or less
lethal strains of the pathogen will predominate. The ratio of infections by more and
less lethal pathogens stabilizes very fast, so that even if measures to prevent further
spread of the disease are put in place as soon as an outbreak is identified (such as
elimination of the animals in an enclosure), the relative frequency of pathogen types
is likely to have changed before action is taken. Other key factors that can also
contribute to the evolution of more lethal strains and their spread in the host
population include transmissibility and mutation rates. Mutation rates are known to
vary greatly between pathogens such as the flu virus, and others such as trematode
worms (liver flukes, shistosomes etc., which reproduce more slowly).

Our models show that an increase in the density of animals on farms or mariculture
facilities will rapidly lead to the dominance of more lethal strains of pathogens if
these can enter the farm and mutations occur to produce more lethal variants.  An
example of this process has recently been described in the mariculture of fish
in~\cite{PuEtAl}. The problem has been recognized in intensive poultry production~
\cite{SlGl}.
and the identification of more virulent trains of Marburg's disease in
chickens~\cite{Whr} may also be an example of this process, and it seems
possible that the advent of very virulent strains of bird flu during the
Spanish Flu epidemic may be linked to high densities of soldiers in
demobilization camps and troop ships, although these men may have been very
susceptible due to their poor condition~\cite{Ma2}. 

Thus the density at which animals are kept should be considered as a risk
factor for the evolution of more lethal diseases. Assuming chaotic character of
movement of animals inside the enclosure where they are kept (as modelled by
independent Brownian motion processes), we discussed how their effective
density affects the key parameter of our model specifying the hosts' encounter
rate and hence eventually determining what pathogen type will predominate.

The ratio of pathogen types will clearly also depend on the transmissibility of
pathogen strains. We have focused on differences in lethality between pathogen
strains, because pathogen strains with increased transmissibility would
obviously become more prevalent. It is interesting to note that our model
showed that, even if the transmissibility of pathogen of one type is lower than
that for the other, the former pathogen can still prevail provided it is
favoured by mutation. Density of the animal hosts may itself increase
transmissibility, due to stress on the hosts caused by
crowding~\cite{Sn}. 
Increased transmissibility might in turn lead to the evolution of increased
lethality~\cite{Ew}. 
An interesting question is whether increased lethality would be linked in most
cases to increased transmissibility. Some previous models~\cite{Ew, Fr, Da}
have assumed that more rapid production of copies of the pathogen
inside the host would both increase the likelihood of transmission to new hosts
and also shorten the life of the host. We did not make any assumptions of that
kind for our model.

Modern animal husbandry often involves a large number of separate enclosures,
each containing a large number of animals at very high densities. Once a
disease is detected in an enclosure, farmers would usually sacrifice or remove
the animals, but pathogens may be carried between enclosures by various
mechanisms, depending on the type of pathogen and the biosecurity practices
followed). A multistage version of our model for this situation suggests that
if pathogens are transferred a number of times, then the evolution of more
lethal pathogens may be very rapid, but the increase in lethality will
eventually lead to the epidemic becoming unsustainable (hosts dying too fast to
be able to transit the pathogen).

We suggest that the outcomes predicted by the mathematical models discussed in
the present paper can carry important messages for animal husbandry, where
there are strong commercial incentives to increase the densities of animals in
enclosures to very high levels, and often very large numbers of enclosures are
built in a single farm.

\medskip

\noindent{\bf Acknowledgments.} K.\,Borovkov was supported by the ARC Centre of
Excellence   for Mathematics and Statistics of Complex Systems.

\end{document}